%% file: instab.tex
\documentclass[10pt,a4paper,leqno]{amsart}

\usepackage{amsmath,amssymb,amsthm}
\usepackage{hyperref}

\def\dis
{\displaystyle}

\def\R{{\mathbb R}}
\def\N{{\mathbb N}}

\def\S{{\mathcal S}}
\def\virgp{\raise 2pt\hbox{,}}
\def\bv{{\bf v}}
\def\bu{{\bf u}}
\def\bw{{\bf w}}
\def\({\left(}
\def\){\right)}
\def\<{\left\langle}
\def\>{\right\rangle}

\def\d{{\partial}}
\def\a{{\tt a}}

\def\g{\gamma}
\def\e{\varepsilon}
\def\l{\lambda}
\def\om{\omega}
\def\si{{\sigma}}

\def\u{{\tt u}}
\def\v{{\tt v}}
\def\w{{\tt w}}

\def\O{\mathcal O}

\theoremstyle{plain}
\newtheorem{theo}{Theorem}[section]
\newtheorem{lem}[theo]{Lemma}
\newtheorem{cor}[theo]{Corollary}
\newtheorem{prop}[theo]{Proposition}
\newtheorem*{hyp}{Assumptions}

\theoremstyle{definition}
\newtheorem{defin}[theo]{Definition}
\newtheorem*{nota}{Notation}

\theoremstyle{remark}
\newtheorem{rema}[theo]{Remark}
\newtheorem*{rema*}{Remark}

\newtheorem{ex}{Example}

\numberwithin{equation}{section}

\begin{document}

\title[Instability for semi-classical Schr\"odinger
equations]{Geometric optics and instability for semi-classical
  Schr\"odinger equations} 
\author[R. Carles]{R{\'e}mi Carles}
\address{Institut math\'ematique de Bordeaux\\33405 
  Talence cedex\\France}
\email{Remi.Carles@math.cnrs.fr}
\thanks{This work was begun while the author was on leave at IRMAR
  (University of Rennes). He wishes to thank this institution for its
  kind hospitality. Support by European network HYKE,
  funded  by the EC as contract HPRN-CT-2002-00282 is also acknowledged.}
\begin{abstract}
We prove some instability
phenomena for semi-classical (linear or) nonlinear Schr\"odinger
equations. For some perturbations of the data, we show that for very
small times, we can neglect the Laplacian, and the 
mechanism is the same as for the corresponding 
ordinary differential equation. Our approach allows smaller
perturbations of the 
data, where the instability occurs for times such that the problem cannot be 
reduced to the study of an o.d.e. 
\end{abstract}
\maketitle

\input{intro}

\input{bkw}
\input{formal}

\input{strong}

\input{caslim}
\input{weak}
\appendix
\input{linear}
\input{ill}

\input{flow}
\bibliographystyle{amsplain}
\bibliography{instab}

\end{document}

%% file: intro.tex
\section{Introduction}
\label{sec:intro}

Consider the semi-classical Schr\"odinger equation:
\begin{equation}\label{eq:schrodcubic}
ih \d_t u^h + \frac{h^2}{2}\Delta u^h =   |u^h|^2u^h
\quad ; \quad u^h(0,x)= a_0(x)\, ,
\end{equation}
where $x\in\R^n$, the parameter $h>0$ goes to zero and the initial
datum $a_0$ is independent of $h$. We prove that
small perturbations of $a_0$ cause divergence of the
corresponding two solutions on small time intervals. For instance,
assume  that $a_0$ is smooth, $a_0\in
\S(\R^n)$. Consider $v^h$ solving \eqref{eq:schrodcubic} with datum $a_0+
h^{1-\frac{1}{N}} a_1$, with $a_1\in \S(\R^n)$, $\rm{Re}(a_0
\overline{a_1})\not \equiv 0$ and $N>0$. Then there exists $c>0$
independent of $h$ such that for $t^h =h^{1/N}$:
\begin{equation*}
\liminf_{h\to 0}\left\| u^h \( t^h\)- v^h \( t^h\)\right\|_{L^2}\ge c\, .
\end{equation*}
Such an instability phenomenon goes in the same spirit as the study
of G.~Lebeau \cite{Lebeau01} (see also \cite{LebeauX00},
\cite{MetivierBourbaki}; see \cite{Lebeau05} for further developments)
for the nonlinear wave 
equation, and followed for instance in  \cite{BGTMRL}, \cite{CCT,CCT2}
(see also the appendix of \cite{BGTENS}) and
\cite{BZ} for nonlinear Schr\"odinger equations. For the above
example, in the case $N<3$, our approach relies on the fact that for
very small 
times, the 
dispersive effects due to the Laplacian are negligible, as in
\cite{Lebeau01} and \cite{CCT2}; a good
approximation to the Schr\"odinger equation is then provided by an ordinary
differential equation, which can be solved explicitly. In the case
$N\ge 3$, the instability mechanism occurs for times such that the
action of the Laplacian is no longer negligible, and the equation
cannot be reduced to an ordinary differential equation. In that case,
our analysis relies on small time properties of solutions of the
compressible Euler 
equation, which describes the semi-classical limit for
\eqref{eq:schrodcubic}.

We also consider weaker nonlinearities in space dimension
$n\ge 2$, with or without
harmonic potential ($\om \ge 0$):
\begin{equation*}
ih \d_t u^h + \frac{h^2}{2}\Delta u^h = \om^2\frac{|x|^2}{2}u^h + h^k
  |u^h|^2u^h\, 
\quad ; \quad u^h(0,x)= a_0(x)\, .
\end{equation*}
In space dimension three,
we can take $k=2$ and $\om>0$, thus recovering the scaling of \cite{BZ}
corresponding to Bose--Einstein condensation in dimension three with repulsive
nonlinearity. Unlike in \cite{BZ} where initial data concentrated at
one point with scale $h$ are considered, we assume that the initial
data  is independent of $h$, $u^h_{\mid t=0}= a_0(x)$. However, the
instability mechanism we describe occurs at a time where the solution
is concentrated. The concentration is due to the presence of the
harmonic oscillator, but the rate of concentration when instability
occurs is smaller than in \cite{BZ} ($h^\alpha$ with $\alpha<1$, see
Section~\ref{sec:weaker} for more details). 

So far, we have considered only cubic nonlinearities. As in
\cite{Grenier98}, we extend the framework to nonlinearities of the
form $f(|u|^2)u$ which are
smooth, repulsive, and cubic at the origin:
\begin{hyp}
  Let $f$ be smooth: $f\in
  C^\infty(\R_+;\R)$, with  $f(0)=0$ and $f'>0$. 
\end{hyp} 
\begin{rema}
The assumption $f(0)=0$ is neutral, since constant potentials for
Schr\"odinger equations can be absorbed by an easy change of unknown
function. 
\end{rema}
\begin{nota}
Let $(\alpha^h)_{0<h\le 1}$ and $(\beta^h)_{0<h\le 1}$ be two families
of positive real numbers. 
\begin{itemize}
\item We write $\alpha^h \ll \beta^h$ if
$\dis \limsup_{h\to 0}\alpha^h/\beta^h =0$.
\item We write $\alpha^h \lesssim \beta^h$ if 
$\dis \limsup_{h\to 0}\alpha^h/\beta^h <\infty$.
\item We write $\alpha^h \approx \beta^h$ if $\alpha^h \lesssim \beta^h$ and
$\beta^h \lesssim \alpha^h$. 
\end{itemize}
\end{nota}
\begin{theo}\label{theo:strong}
Let $n\ge 1$, $a_0 ,\widetilde a_0^h\in\S(\R^n)$,  $\phi_0\in
C^\infty(\R^n;\R)$, where $a_0$ and $\phi_0$ are 
independent of $h$, and $\nabla \phi_0 \in H^s(\R^n)$ for every $s\ge
0$.  For $\omega \ge 0$, let $u^h$ and 
$v^h$ solve the initial value problems:
\begin{align*}
ih \d_t u^h + \frac{h^2}{2}\Delta u^h &= \omega^2\frac{|x|^2}{2}u^h +
f\(|u^h|^2\)u^h 
\quad ; \quad u^h(0,x)= a_0(x)e^{i\phi_0(x)/h}\, .\\
ih \d_t v^h + \frac{h^2}{2}\Delta v^h &= \omega^2\frac{|x|^2}{2}v^h +
f\(|v^h|^2\)v^h 
\quad ; \quad v^h(0,x)= \widetilde a_0^h(x)e^{i\phi_0(x)/h}\, .
\end{align*}
Assume that there exists $N\in \N$ and $h^{1-\frac{1}{N}}\ll \delta^h
\ll 1$ such that:
\begin{equation}\label{eq:hyppert}
\begin{aligned}
  \left\|a_0 -\widetilde
    a_0^h\right\|_{H^s} \approx \delta^h\, ,\quad
  \forall s\ge 0\quad ;\quad
\limsup_{h\to 0}\left\| \frac{\operatorname{Re}(a_0-\widetilde
    a_0^h)\overline{a_0}}{\delta^h}  \right\|_{L^\infty(\R^n)} \not =0.
\end{aligned}
\end{equation}
Then we can find $0<t^h\ll 1$ such that:
$\displaystyle \left\| u^h(t^h) - v^h(t^h) \right\|_{L^2}\gtrsim
1$. More precisely, this mechanism occurs as soon as $t^h \delta^h
\gtrsim h$. 
In particular, 
\begin{equation*}
\frac{\left\| u^h - v^h
  \right\|_{L^\infty([0,t^h];L^2)}}{\left\| u^h_{\mid t=0} - v^h_{\mid
  t=0}
  \right\|_{L^2}}\to +\infty \quad \text{as }h\to 0\, .
\end{equation*}
\end{theo}
\begin{rema*}
  We state assumptions in $H^s$ for
  \emph{every} $s\ge 0$. It will 
  appear in the proof that choosing $s$ large enough would
  suffice. Similarly, the assumption $a_0, \widetilde a_0^h\in
  \S(\R^n)$ is not necessary in our proof. 
\end{rema*}
\begin{rema*}
  The second part of the assumption \eqref{eq:hyppert} can be viewed
  as a polarization condition. We could remove it with essentially the
  same proof as below, up to demanding $h^{1/2-1/N}\ll \delta^h \ll 1$. 
\end{rema*}
The above result can be applied in the following cases:
\begin{ex}
  Consider $a_0,b_0\in \S(\R^n)$ independent of $h$, such that
  $\operatorname{Re}(\overline{a_0}b_0)\not \equiv 0$, and take $\widetilde
  a_0^h = a_0  +\delta^h b_0$. 
\end{ex}
\begin{ex}
  Consider $a_0\in \S(\R^n)$ independent of $h$ and $x^h\in \R^n$. We
  can take $\widetilde  a_0^h(x) = a_0 (x-x^h)$, provided that 
$|x^h|=\delta^h$ and
\begin{equation*}
  \limsup_{h\to 0} \left\| \frac{x^h}{|x^h|}\cdot \nabla\(
    |a_0|^2\)\right\|_{L^\infty} \not =0. 
\end{equation*}
\end{ex}
The above result addresses perturbations which satisfy in particular
$\delta^h \gg 
h$. This excludes the standard WKB data of the form
\begin{equation*}
u^h\big|_{t=0}= a_h(x)e^{i\phi_0(x)/h}\, ,\text{ where }a_h \sim
a_0+ha_1+h^2a_2+\ldots 
\end{equation*}
In that case, a perturbation of $a_1$ is relevant at
time $t^h\approx 1$, and the previous result is essentially sharp:
\begin{prop}\label{prop:caslim}
 Let $n\ge 1$, $a_0,a_1 \in\S(\R^n)$,
 $\phi_0 \in C^\infty(\R^n;\R)$ independent of $h$, with $\nabla
 \phi_0 \in H^s(\R^n)$ for every $s\ge 
0$. Assume that
 $\operatorname{Re}( a_0\overline{a_1})\not \equiv 0$.   For $\omega
 \ge 0$, let $u^h$ and
$v^h$ solve the initial value problems:
\begin{align*}
ih \d_t u^h + \frac{h^2}{2}\Delta u^h &= \omega^2\frac{|x|^2}{2}u^h +
f\(|u^h|^2\)u^h 
\quad ; \quad u^h(0,x)= a_0(x)e^{i\frac{\phi_0(x)}{h}}\, .\\
ih \d_t v^h + \frac{h^2}{2}\Delta v^h &= \omega^2\frac{|x|^2}{2}v^h +
f\(|v^h|^2\)v^h 
\quad ; \quad v^h(0,x)= \(a_0(x)+ha_1(x)\)e^{i\frac{\phi_0(x)}{h}}.
\end{align*}
Then  for any $\tau^h \ll 1$, $\displaystyle \left\| u^h - v^h
\right\|_{L^\infty([0,\tau^h];L^2)}\ll 1 $, and for
$t>0$ independent of $h$ arbitrarily small: 
$\displaystyle \left\| u^h(t) - v^h(t) \right\|_{L^2}\gtrsim 1.$
\end{prop}

For weaker nonlinearities, we have the following result. The
notation $\e$ for the small parameter instead of $h$ is neither a 
mistake nor a coincidence (see Section~\ref{sec:weaker}). 

\begin{cor}\label{cor:weak}
Let $n\ge 2$, $1<k<n$, $a_0,\widetilde a_0^\e\in\S(\R^n)$   and $\om
\ge 0$. Let $u^\e$ and  
 $v^\e$ solve the initial value problems:
\begin{align*}
i\e \d_t u^\e + \frac{\e^2}{2}\Delta u^\e &= \om^2\frac{|x|^2}{2}u^\e +
f\(\e^k|u^\e|^2\)u^\e\, 
\ ; \ u^\e\big|_{t=0}= a_0(x)e^{i\phi_0(x)/\e}\, ,\\
i\e \d_t v^\e + \frac{\e^2}{2}\Delta v^\e &= \om^2\frac{|x|^2}{2}v^\e +
f\(\e^k|v^\e|^2\)v^\e\, 
\ ; \ v^\e\big|_{t=0}= \widetilde a_0^\e(x)e^{i\phi_0(x)/\e}\, .
\end{align*}
Assume that:
\begin{itemize}
\item Either:
there exists $N\in \N$ and $\e^{1-\frac{k}{n}-\frac{1}{N}}\ll \delta^\e
\ll 1$ such that:
\begin{align*}
  \left\|a_0 -\widetilde
    a_0^\e\right\|_{H^s} \approx \delta^\e\, ,\quad
  \forall s>0\quad ;\quad  
\limsup_{\e\to 0}\left\| \frac{\operatorname{Re}(a_0-\widetilde
    a_0^\e)\overline{a_0}}{\delta^\e}  \right\|_{L^\infty(\R^n)} \not =0.
\end{align*}
\item Or: $\widetilde a_0^\e = a_0 +\e^{1-\frac{k}{n}}a_1$, with
  $a_1\in \S(\R^n)$,  $\operatorname{Re}( a_0\overline{a_1})\not
  \equiv 0$. 
\end{itemize}
$1.$ $\om =0$: let $\phi_0(x)=-\frac{|x|^2}{2}$. There
exist $T^\e \to 1^{-}$  and $0<\tau^\e\ll 1$ such that:
\begin{equation}\label{eq:instabweak}
\left\| u^\e - v^\e \right\|_{L^\infty([0,T^\e];L^2)}\ll 1\quad
;\quad 
\left\| u^\e - v^\e \right\|_{L^\infty([0,T^\e+\tau^\e];L^2)}\gtrsim 1\, .
\end{equation}
$2.$ $\om >0$: let $\phi_0\equiv 0$. There
exist $T^\e\to \frac{\pi}{2\om}^{-}$ and $0<\tau^\e\ll 1$ such that
\eqref{eq:instabweak} holds. 
\end{cor}
\begin{ex}
If $n=3$, $k=2$ and the nonlinearity is cubic, we consider:
\begin{equation*}
i\e \d_t u^\e + \frac{\e^2}{2}\Delta u^\e = \om^2\frac{|x|^2}{2}u^\e +
\e^2|u^\e|^2 u^\e\,. 
\end{equation*}
Then perturbations of order $\delta^\e$ with $1\gg \delta^\e \gg
\e^{1/3-1/N}$ cause instability. On the other hand, this
phenomenon does not occur for the same equation in space
dimension two, and there is stability for a large class of
initial data (see \cite{CaIHP,CK}). 
\end{ex}
\begin{rema}
We have $T^\e\to 1^{-}$ in the first case because of the initial
quadratic oscillations. In the linear case, such oscillations cause
focusing at the origin at time $t=1$ (see e.g. \cite{CaIUMJ}). We will
see that the instability mechanism occurs when the solution is no
longer of order $\O(1)$ and is already concentrated at scale
$1-T^\e$. In the case $\om >0$, a similar phenomenon occurs
without initial phase because the action of the harmonic oscillator is
similar. In both cases, taking $\phi_0(x)=-b|x|^2$ and modulating $b$,
we could have the instability mechanism occur near any time $T>0$, and
not only $1$ or $\frac{\pi}{2\om}$. 
\end{rema}
\begin{rema} 
In Corollary~\ref{cor:weak}, the assumption $k<n$ is crucial. When $k=n$,
the above result is no longer true (see \cite{CaIUMJ,CaIHP,CFG,CK} for an
homogeneous 
nonlinearity, with or without harmonic potential). From the point of
view of geometrical optics, assuming  $k<n$ amounts to considering a
super-critical r\'egime if a caustic reduced to a point appears. This
goes in the 
spirit of the formal computations of \cite{HK87}, and of the papers
\cite{JMR95,JMR00a,JMR00b}, \cite{Lebeau01,MetivierBourbaki},
\cite{CaIUMJ,Ca-p,CFG}.  
\end{rema}
The above results show in particular that computing directly 
semi-classical limits of nonlinear Schr\"odinger equations by numerical
methods is highly challenging. \\

Compare with other results on instability. In
\cite{MetivierBourbaki,CCT2,BGTENS,BZ}, the perturbation are of order
$|\ln h|^{-\theta}$ for some $\theta >0$. Then instability occurs at
time of order $h |\ln h|^{\theta '}$. Our analysis allows smaller
perturbations, and the instability occurs a little later. \\

The paper is organized as follows.  Section~\ref{sec:bkw} is devoted
to a general discussion on WKB methods for nonlinear Schr\"odinger
equations. In Section~\ref{sec:formal}, we give heurisitc arguments to
prepare the proof of Theorem~\ref{theo:strong} and
Corollary~\ref{cor:weak}. The proof of Theorem~\ref{theo:strong} is
completed in Section~\ref{sec:strong}, and
Proposition~\ref{prop:caslim} is established in
Section~\ref{sec:caslim}. Corollary~\ref{cor:weak} is shown in
Section~\ref{sec:weaker}. In a first appendix, we exhibit some notion of
instability in the linear case; in a second appendix, we show how some
results of \cite{CCT,CCT2} can be recovered from semi-classical
analysis. Finally in Appendix~\ref{sec:flow}, we establish
the following result:
\begin{cor}\label{cor:flow}
   Let $n\ge 3$. Consider the cubic,
  defocusing Schr\"odinger equation:
  \begin{equation}
    \label{eq:NLScubic}
    i\d_t u +\frac{1}{2}\Delta u = |u|^2 u\quad ; \quad u_{\mid
    t=0}=u_0\, .
  \end{equation}
Denote $s_c = \frac{n}{2}-1$. Let $0<s<s_c$. We can find a
  family $(u_0^\e)_{0<\e \le 
  1}$ in $\S(\R^n)$ with 
\begin{equation*}
  \|u_0^\e\|_{H^s(\R^n)} \to 0 \text{ as }\e \to 0\, ,
\end{equation*}
and $0<t^\e \ll 1$ such that the solution $u^\e$ to
\eqref{eq:NLScubic} associated to $u_0^\e$ satisfies: 
\begin{equation*}
  \|u^\e(t^\e)\|_{H^{k}(\R^n)} \to +\infty \text{ as }\e \to 0\, , \ \forall
 k\in \left]\frac{s}{\frac{n}{2}-s}, s\right]\, .
\end{equation*}
In particular, for any $t>0$, the map $u_0 \mapsto u(t)$ given by
\eqref{eq:NLScubic} fails to be continuous at $0$ from $H^s$ to
$H^{k}$. 
\end{cor}
\begin{rema}
This result can be viewed as a weak version of the analog result to
\cite{Lebeau05} for the cubic, defocusing Schr\"odinger equation. An
important difference though is that we consider a \emph{sequence} of
initial data, while in \cite{Lebeau05}, G.~Lebeau considers the weak
solution associated to a \emph{fixed} initial data. What prevents us from
filling this gap is the finite speed of propagation which is used in
\cite{Lebeau05} for the wave equation, and is not available for
Schr\"odinger equations.   
\end{rema}

%% file: bkw.tex
\section{WKB methods for nonlinear Schr\"odinger equations}
\label{sec:bkw}

Consider the initial value problem, for $x\in \R^n$:
\begin{equation}
  \label{eq:nlssemi}
  ih \d_t u^h +\frac{h^2}{2}\Delta u^h = h^\kappa |u^h|^2
  u^h\quad ; \quad u^h_{\mid t=0} = a_0^h(x)e^{i\phi_0(x)/h}\, . 
\end{equation}
The aim of WKB methods is to describe $u^h$ in the limit $h \to 0$,
when $\phi_0$ does not depend on $h$, and 
$a_0^h$ has an asymptotic expansion of the form:
\begin{equation*}
  a_0^h (x)\sim a_0(x) +h a_1(x)+h^2 a_2(x) +\ldots
\end{equation*}
The parameter $\kappa \ge 0$ describes the strength of a coupling
constant, which makes nonlinear effects more or less important in the
limit $h \to 0$; the larger $\kappa$, the weaker the nonlinear
interactions. Note that since we consider an homogeneous nonlinearity,
this amounts to considering the case where the coupling constant is
$1$, with initial data of order $h^{\kappa /2}$.

An interesting feature of \eqref{eq:nlssemi} is that one does not
expect the creation of harmonics. The WKB methods consist in seeking
an approximate solution to 
\eqref{eq:nlssemi} of the form:
\begin{equation}\label{eq:defBKW}
  u^h(t,x)\sim \(\a^{(0)}(t,x) + h \a^{(1)}(t,x) + h^2
  \a^{(2)}(t,x)+\ldots \) e^{i\phi(t,x)/h}\, .
\end{equation}
For such an expansion to be available with profiles $\a^{(j)}$
independent of $h$, it is reasonable to assume that $\kappa$ is an
integer, $\kappa \in \N$. 
One must not expect this approach to be valid when caustics are
formed: roughly speaking, when a caustic appears, all the terms
$\phi$, $\a^{(0)}$, $\a^{(1)}$, \ldots become singular. In this paper,
we always consider times preceding this break-up.

\subsection{Notion of criticality}
\label{sec:crit}

If $\kappa \ge 2$, then nonlinear effects are negligible at
leading order in WKB methods. On the other hand, if $\kappa =1$
(weakly nonlinear geometric 
optics), then nonlinear effect are relevant at leading order. The
present discussion is 
formal, its aim being to prepare the study of the case $\kappa =0$.  
 
When $\kappa \ge 2$, plugging the asymptotic expansion \eqref{eq:defBKW}
into \eqref{eq:nlssemi} yields formally:
\begin{equation*}
  \begin{aligned}
    \d_t \phi + \frac{1}{2}|\nabla \phi|^2 =0 \quad & ;\quad
    \phi_{\mid t=0} = \phi_0 \, .\\
\d_t \a^{(0)} + \nabla \phi \cdot \nabla \a^{(0)}
    +\frac{1}{2}\a^{(0)}\Delta \phi =0 \quad & ;\quad
    \a^{(0)}_{\mid t=0} = a_0\, .
  \end{aligned}
\end{equation*}
The first equation is the well known eikonal equation, which describes
the geometry of the propagation. If $\phi_0$ is smooth, it has a
smooth solution, locally in time. This solution may become singular in
finite time, this phenomenon being the formation of a caustic. 

The second equation is a transport
equation, which is simply an ordinary differential equation for  the
leading order amplitude  along the rays of geometrical optics. To see
this, introduce a parametrization of these rays:
\begin{equation*}
  \frac{d}{dt}X_t(x) = \nabla \phi\(t,X_t(x)\) \quad ;\quad X_0(x)=x
  \, ,
\end{equation*}
and the Jacobi determinant:
$\displaystyle  J_t(x) =\operatorname{det}\nabla X_t(x)$.
It is well defined and smooth so long as no caustic appears.
The break-up time $t_c>0$, if any, is such that there exists $x_c$ such
that $J_{t_c}(x_c)=0$. 
The transport equation for $\a^{(0)}$ is  the trivial
ordinary differential equation:
\begin{equation*}
  \frac{d}{dt} \( \a^{(0)}\( t,X_t(x)\)\sqrt{J_t(x)}\)=0\, .
\end{equation*}
From this, we easily see that when a caustic appears, not only $\phi$
becomes singular, but also $\a^{(0)}$, since $J_t(x)$ goes to zero at
the caustic. 


The value $\kappa =1$ is critical as far as leading order phenomena
are concerned: the transport equation for $\a^{(0)}$ is then nonlinear,
\begin{equation}\label{eq:weak}
  \begin{aligned}
    \d_t \phi + \frac{1}{2}|\nabla \phi|^2 &=0 \quad  ;\quad
    \phi_{\mid t=0} = \phi_0 \, .\\
\d_t \a^{(0)} + \nabla \phi \cdot \nabla \a^{(0)}
    +\frac{1}{2}\a^{(0)}\Delta \phi &=-i \left|\a^{(0)}\right|^2 \a^{(0)}
    \quad  ;\quad 
    \a^{(0)}_{\mid t=0} = a_0\, .
  \end{aligned}
\end{equation}
On the other hand, the eikonal equation is still the same as in the
linear case, hence the term ``weakly nonlinear'' (see also
\cite{RauchUtah} and references therein).  The correctors 
$(\a^{(j)})_{j\ge 2}$ solve  linear transport
equations. With the above notations, the
nonlinear transport equation is again an ordinary differential
equation along rays:
\begin{equation*}
  \frac{d}{dt} \( \a^{(0)}\( t,X_t(x)\)\sqrt{J_t(x)}\)=-i
\left|\a^{(0)}\( t,X_t(x)\)\right|^2
\a^{(0)}\( t,X_t(x)\)\sqrt{J_t(x)}\, .
\end{equation*}
This ordinary differential
equation is of the form $\dot y = iV y$, where the nonlinear potential
$V$ is real-valued. In particular, the modulus of $y$ is constant, and
we just have to solve a \emph{linear} differential
equation. Thus, leading order
nonlinear effects are measured by a (nonlinear) phase shift, which may
be compared to the phenomenon of \emph{phase self-modulation} in laser
physics (see e.g. \cite{ZS,Boyd,Donnat}).

\subsection{Super-critical case}
\label{sec:super}

In the super-critical case $\kappa =0$, the nonlinearity is
present in the eikonal equation: the hierarchy of the case $\kappa=1$
is shifted, so that the corrector $\a^{(1)}$ is present in the
transport equation for $\a^{(0)}$. As noted in
\cite{PGX93},  the system for the phase $\phi$ and the
amplitudes 
$\a^{(0)}$, $\a^{(1)},\ldots$ is not closed (see also
\cite{CheverryBullSMF,CGM04,CG05}). For instance, we find:
\begin{equation}\label{eq:strong}
  \begin{aligned}
    \d_t \phi + \frac{1}{2}|\nabla \phi|^2 +\left|\a^{(0)}\right|^2
    &=0 \, .\\
\d_t \a^{(0)} + \nabla \phi \cdot \nabla \a^{(0)}
    +\frac{1}{2}\a^{(0)}\Delta \phi &=-2i
    \a^{(0)}\operatorname{Re}\(\a^{(0)}\overline{\a^{(1)}} \)\, .\\
\d_t \a^{(1)} + \nabla \phi \cdot \nabla \a^{(1)}
    +\frac{1}{2}\a^{(1)}\Delta \phi &=\frac{i}{2}\Delta
    \a^{(0)}-i\left|\a^{(1)}\right|^2 a^{(0)}\\
 -2i
    \a^{(1)}&\operatorname{Re}\(\a^{(0)}\overline{\a^{(1)}} \)
-2i
    \a^{(0)}\operatorname{Re}\(\a^{(0)}\overline{\a^{(2)}} \) 
   \, .
  \end{aligned}
\end{equation}
However, as pointed out in \cite{PGX93}, the phase $\phi$ can be found when
considering: 
$$(\rho ,v)= \(|\a^{(0)}|^2, \nabla
\phi\).$$ 
Indeed, it solves the compressible, isentropic Euler
equation:
\begin{equation}\label{eq:euler}
  \begin{aligned}
    \d_t \rho +\operatorname{div}\(\rho v\) &= 0\quad ;\quad \rho_{\mid
    t=0}=|a_0|^2\, .\\
\d_t v +v\cdot \nabla v + \nabla \rho &= 0\quad ;\quad v_{\mid
    t=0}=\nabla\phi_0\, .
  \end{aligned}
\end{equation}
For smooth initial data decaying to zero at infinity, this system as a
smooth solution locally in time \cite{Majda,MUK86,JYC90}. In general, 
finite time blowup occurs \cite{MUK86,JYC90,Xin}, but not always
\cite{Grassin}; the known results depend on the propagation of 
the initial velocity by the (multi-dimensional) Burgers' equation. 

Once $(\rho,v)$ is
determined, $\d_t \phi$ is given by the eikonal equation; this yields
$\phi$. Note that knowing $(\rho,v)$ suffices to compute important
quadratic quantities such as Wigner measures.
To complete the closure of the system, and provided that the
leading order amplitude $a_0$ is nowhere zero, one may consider a
generalized Madelung transform (see \cite{PGX93}), which we do not
describe here.

\subsection{Justification on small time intervals}
\label{sec:justifsmall}
 
Justifying geometric optics in the super-critical case is, in general,
an open problem. However, as noticed in \cite{PGX93} and exploited in
many other works (see
e.g. \cite{Kuksin95,MetivierBourbaki,CCT2,BGTENS,BZ}), if one studies
this limit on time 
intervals of the form $[0, c_0 h|\ln h|]$ for some $c_0>0$, then the
problem is simpler. Consider the more general nonlinear Schr\"odinger
equation in 
$\R^n$: 
\begin{equation}
  \label{eq:NLSgensemi}
  ih\d_t u^h +\frac{h^2}{2}\Delta u^h = \omega |u^h|^{2\si}u^h \quad ;
\quad u^h(0,x)=a_0(x)\, ,
\end{equation}
where $\omega \in \R\setminus \{0\}$ and $\si\in \N\setminus\{0\}$. We
consider the case $\phi_0\equiv 0$ to prove that in this case, one can
choose an 
approximate which is even simpler than the one given by
\eqref{eq:euler}. 
Formally, $u^h$ is formally approximated by $a
e^{i\phi /h}$ where:
\begin{equation*}
\begin{aligned}
\d_t \phi +\frac{1}{2}|\nabla \phi|^2 + \omega |a|^{2\si} =0 \quad ;\quad
\phi_{\mid t=0} =0\, .\\
\d_t a +\nabla \phi\cdot \nabla a +\frac{1}{2}a\Delta \phi =0\quad
; \quad a_{\mid t=0}=a_0\, . 
\end{aligned}
\end{equation*}
Looking at Taylor expansions for $\phi$ and $a$ as $t\to 0$, we see that
\begin{equation*}
  a(t,x) = a_0(x) + \O(t^2)\quad ;\quad \phi(t,x) = -t \omega
  |a_0(x)|^{2\si}+\O(t^3)\, .
\end{equation*}
We prove that $a(t,x)
e^{i\phi(t,x) /h}$ can approximated by $a_0(x) e^{-it \omega
  |a_0(x)|^{2\si}/h}$ on some time interval of the form $[0, c_0 h|\ln
h|]$. Call $\varphi^h$ the latter function. It solves the ordinary
differential equation: 
\begin{equation}\label{eq:varphi}
  ih\d_t \varphi^h  = \omega
 |\varphi^h|^{2\si}\varphi^h \quad ;
\quad \varphi^h_{\mid t=0}=a_0\, .
\end{equation}
We prove that if $c_0$ is sufficiently small, then $\varphi^h$ is a
good approximation of $u^h$ on $[0, c_0 h|\ln
h|]$. Let $w^h = u^h-\varphi^h $. It solves:
\begin{equation}
  \label{eq:1435}
  ih\d_t w^h +\frac{h^2}{2}\Delta w^h= \omega\( F(w^h
  +\varphi^h) -  F(\varphi^h)\) +\O(h^2)+\O(t^2)+\O(h t)\, ,
\end{equation}
with $ w^h_{\mid t=0}=0$, where we have set $F(z)=|z|^{2\si}z$. The
$\O(h^2)$ term corresponds to the fact that we consider only the first
two terms of a WKB analysis, the term $\O(t^2)$ stems from the
approximation of the phase for small times, and $\O(ht)$ from the
approximation of the amplitude for small times. To be more precise, we
must say that these source terms are measured in $L^2\cap L^\infty
(\R^n)$. When measured in $H^k$, they must be multiplied by a factor
of order $ 1 + (t/h)^k$, due to the differentiation of the phase. 
For $k\ge 0$, we have:
\begin{align*}
  \|w^h\|_{L^\infty([0,t];H^k)} \lesssim & \frac{1}{h}\left\|F(w^h
  +\varphi^h) -  F(\varphi^h)  \right\|_{L^1([0,t];H^k)} \\
&+
  \frac{1}{h}\int_0^t \(h^2  + 
  s^2 +h s\)\left\langle\frac{s}{h}\right\rangle^kds .
\end{align*}
At least for $\si$ integer, we have, when $k>n/2$:
\begin{align*}
 \left\|F(w^h(t)
  +\varphi^h(t)) -  F(\varphi^h(t))  \right\|_{H^k}&\lesssim \(
  \|w^h(t)\|^{2\si}_{H^k} +
  \|\varphi^h(t)\|^{2\si}_{H^k}\)\|w^h(t)\|_{H^k}\\
&\lesssim \(
  \|w^h(t)\|^{2\si}_{H^k} +
  \left\langle\frac{t}{h}\right\rangle^{2\si k}\)\|w^h(t)\|_{H^k}.
\end{align*}
On any time interval where we have, say, $\|w^h\|_{H^k}\le 1$, we infer:
\begin{equation*}
  \|w^h\|_{L^\infty([0,t];H^k)} \le
  \frac{C}{h}\int_0^t \left\langle\frac{s}{h}\right\rangle^{2\si k} 
  \| w^h (s)\|_{H^k}ds + C_1\int_0^t \(h  +
  \frac{s^2}{h} +s\)\left\langle\frac{s}{h}\right\rangle^k ds.
\end{equation*}
Gronwall lemma yields:
\begin{align*}
  \|w^h\|_{L^\infty([0,t];H^k)} &\lesssim \int_0^t \(h  +
  \frac{s^2}{h}
  +s\)\left\langle\frac{s}{h}\right\rangle^k\exp\(\frac{C}{h}\int_s^t
  \left\langle\frac{\tau}{h}\right\rangle^{2\si k}  
  d\tau\)ds.
\end{align*}
Let $t^h = c_0 h |\ln h|^\theta$:
\begin{align*}
  \|w^h\|_{L^\infty([0,t^h];H^k)} &\lesssim \exp\( Cc_0 |\ln h|^\theta
  \left\langle \ln h \right\rangle^{2\si k\theta}  \)\int_0^{t^h} \(h s +
  \frac{s^3}{h}
  +s^2\)\left\langle\frac{s}{h}\right\rangle^k ds\\
&\lesssim \exp\( Cc_0 |\ln h|^\theta
  \left\langle c_0\ln h \right\rangle^{2\si k\theta}  \)h^2 |\ln h|^{4\theta}.
\end{align*}
For $\theta =(1+2\si k)^{-1}$ and $c_0$ sufficiently small, this
yields:
\begin{equation*}
  \|w^h\|_{L^\infty([0,t^h];H^k)} \lesssim h |\ln h|^{\frac{4}{1+2\si k}}.
\end{equation*}
We can then conclude with a continuity argument, for $h$
sufficiently small:
\begin{prop}\label{prop:bkwlog}
  Let $n\ge 1$, $\omega \in \R\setminus \{0\}$, $\si>0$ an integer,
  and $a_0\in \S(\R^n)$. Fix $k>n/2$. Then we can find $c_0,c_1,\theta>0$
  independent of $h\in 
  ]0,1]$ such that $u^h$ and $\varphi^h$, solutions to
  \eqref{eq:NLSgensemi} and \eqref{eq:varphi} respectively, satisfy:
\begin{equation*}
 \|u^h-\varphi^h\|_{L^\infty([0,c_0 h|\ln h|^\theta];H^k)}  \lesssim h
  |\ln h|^{c_1}\, . 
\end{equation*}
\end{prop}

%% file: formal.tex
\section{Instability: formal computations}
\label{sec:formal}

In this section, we show how to reduce the proof of
Theorem~\ref{theo:strong} to the justification of super-critical
nonlinear geometric optics on a time interval which is independent of
$h$. This formal approach would remain valid for a larger class of
nonlinearities, not necessarily defocusing and cubic at the origin. 

\subsection{The o.d.e. mechanism}
\label{sec:ode}

Consider the general Schr\"odinger equation with data independent of $h$:
\begin{equation}\label{eq:schrodgen2}
ih \d_t u^h + \frac{h^2}{2}\Delta u^h =   f\(|u^h|^2\)u^h
\quad ; \quad u^h(0,x)= a_0(x)\, .
\end{equation}
The instability mechanism we sketch in this section is valid for initial data
which are not highly oscillatory: $\phi_0\equiv 0$. We study a more
general framework, corresponding to Theorem~\ref{theo:strong} in
Section~\ref{sec:nonode} below. 

Following WKB methods, seek $u^h$ such that $u^h \sim a e^{i\phi /h}$
as $h\to 0$. Plugging this ansatz into \eqref{eq:schrodgen2} and
canceling $\O(1)$ and $\O(h)$ terms yields:
\begin{equation}\label{eq:bkwlim}
\begin{aligned}
&\d_t \phi +\frac{1}{2}|\nabla \phi|^2 +f\(|a|^2\)=
  0\quad ;\quad 
  \phi(0,x)=0\,. \\
& \d_t a + \nabla \phi\cdot \nabla a +\frac{1}{2}a \Delta \phi
=0 \quad ; \quad a(0,x)=a_0(x)\, .
\end{aligned}
\end{equation}
As $t\to 0$, approximate $\phi$ and $a$ by their Taylor expansion:
\begin{equation*}
\phi(t,x)\sim \sum_{j\ge 1}t^{2j-1}\phi_j(x)\quad ; \quad a(t,x)\sim
\sum_{j\ge 0}t^{2j}{a}_j(x) \, .
\end{equation*}
Note that for $j=0$, the notations are consistent. The fact that only
odd (resp. even) powers of $t$ appear in the expansion for $\phi$
(resp. $a$) is due to the assumption $\phi_{\mid t=0}\equiv
0$. Plugging these 
formal series into \eqref{eq:bkwlim}, we find:
\begin{align*}
\phi_1  = -f\(|a_0|^2\)\quad ; \quad
2a_1 =- \nabla 
\phi_1\cdot \nabla a_0 -\frac{1}{2}a_0 \Delta \phi_1\, .
\end{align*}
Thus, $a_1$ is the first term where the presence of the
Laplacian becomes relevant: let $\u_1^h(t,x)= a_0(x)\exp\(
it\phi_1(x)/h\)$. It solves the ordinary differential
equation:  
\begin{equation}\label{eq:edo}
ih \d_t \u_1^h  =   f\(|\u_1^h|^2\)\u_1^h
\quad ; \quad \u_1^h(0,x)= a_0(x)\, ,
\end{equation}
where $x$ is now just a parameter. Assume that for some time interval 
$[0,T^h]$, WKB method provides a good approximation for $u^h$ in
$L^2(\R^n)$: 
\begin{equation*}
\left\| u^h - a e^{i\phi/h}\right\|_{L^\infty([0,T^h];L^2)}\to 0 \text{ as }
h\to 0\, . 
\end{equation*}
On the other hand, we can approximate $a e^{i\phi/h}$ by $\u_1^h$ if 
\begin{equation*}
\left\|a e^{i\phi/h} - a_0e^{it\phi_1/h}
\right\|_{L^\infty([0,T^h];L^2)}\to 0 \text{ as } 
h\to 0\, .
\end{equation*}
If $T^h \to 0$, which we may assume in view of Th.~\ref{theo:strong},
then approximating $a$ by $a_0$ is not a 
problem. We have to be more careful with the phase, because of the
division by $h$. Formally, the above limit holds if
\begin{equation*}
\exp\(i t^3 /h\)\to 1 \text{ as } 
h\to 0\, .
\end{equation*}
If $t^h\le T^h$ is such that $(t^h)^3\ll h$, then we expect:
\begin{equation}\label{eq:redu}
\left\| u^h - \u_1^h\right\|_{L^\infty([0,t^h];L^2)}\to 0 \text{ as }
h\to 0\, .
\end{equation} 
Now let $v^h$ solve \eqref{eq:schrodgen2} with initial data $v^h(0,x)=
\widetilde a_0^h (x)$, where $\widetilde a_0^h$ satisfies
\eqref{eq:hyppert}; the
assumption on $\delta^h$  
will appear later. Let $\v_1^h$ be the
solution of the corresponding ordinary differential
equation. Similarly, we expect:
\begin{equation}\label{eq:redv}
\left\| v^h - \v_1^h\right\|_{L^\infty([0,t^h];L^2)}\to 0 \text{ as }
h\to 0\, .
\end{equation}
An instability like in Th.~\ref{theo:strong} then stems from an
instability at the o.d.e. level:
\begin{equation*}
\v_1^h(t,x) = \widetilde a_0^h(x)\exp\(i t\widetilde
\phi_1^h(x)/h \)\, ,
\end{equation*}
where $\widetilde
\phi_1^h = -f\(|\widetilde a_0^h|^2\)$. We have obviously
\begin{equation*}
\left\| \v_1^h - a_0\exp\(i t\widetilde
\phi_1^h/h \) \right\|_{L^\infty([0,t^h];L^2)}\to 0 \text{ as }
h\to 0\, ,
\end{equation*} 
as soon as $\delta^h\ll 1$. Instability comes from the
phase:
\begin{equation*}
\u_1^h(t,x)-\v_1^h(t,x) \sim a_0(x) \( e^{it\phi_1(x)/h} -
e^{it\widetilde \phi_1^h(x)/h} \)\, .
\end{equation*}
Using Taylor formula for $f$, we have:
\begin{equation*}
\widetilde
\phi_1^h =  - f\(|a_0|^2\)-2\delta^h\operatorname{Re}\( (a_0 -\widetilde
a_0^h)\overline {a_0}\)f'\(|a_0|^2\)
+\O\(\(\delta^h\)^2\)\,.
\end{equation*}
Since $f'>0$, we infer from \eqref{eq:hyppert}:
\begin{equation*}
\widetilde
\phi_1^h(x) = \phi_1(x)+\delta^hc(x)
+o\(\delta^h\)\,,
\end{equation*}
where the function $c$ does not depend on $h$ and is not identically
zero on the support of $a_0$. For $t\delta^h \approx h$, we infer:
\begin{equation*}
\left|\u_1^h(t,x)-\v_1^h(t,x)\right| \sim |a_0(x)| \left| 
e^{it\delta^hc(x) /h} -1\right|\, .
\end{equation*}
This has a nonzero limit as $h\to 0$ since $t\delta^h \approx h$. 
The only constraint we imposed so far was $t^h \ll  h^{1/3}$, so taking
$h^{2/3} \ll \delta^h \ll 1$ predicts an instability as stated in
Th.~\ref{theo:strong}. To prove Th.~\ref{theo:strong}, we must
establish \eqref{eq:redu} and \eqref{eq:redv} for 
suitable $t^h$.

\begin{rema}\label{rem:proj}
In view of \cite{BZ}, introduce the complex projective distance:
\begin{equation*}
u_j\in L^2(\R^n),\quad d_{\rm pr}(u_1,u_2):= \arccos \(
\frac{|\<u_1,u_2\>|}{\|u_1\|_{L^2} \|u_2\|_{L^2}}\)\,.
\end{equation*}
Then we can check that up to demanding $t^h\ll h^{1/3}$ and $t^h
\delta^h \gg h$ (these conditions can be satisfied for $
h^{2/3}\ll\delta^h\ll 1$),
and provided that \eqref{eq:redu} and 
\eqref{eq:redv} hold: 
\begin{equation*}
\frac{d_{\rm pr} \(u^h(t^h),v^h(t^h)\)}{d_{\rm pr}
  \(u^h(0),v^h(0)\)}\approx
\frac{d_{\rm pr} \(\u^h(t^h),\v^h(t^h)\)}{d_{\rm pr}
  \(\u^h(0),\v^h(0)\)}\to +\infty\quad \text{as }h\to 0\, .
\end{equation*}
\end{rema}
\begin{rema}
  We prove in Appendix~\ref{sec:linear} a result in a similar spirit
  for linear equations. 
\end{rema}
\begin{rema}
  We show in Appendix~\ref{sec:cct} how this analysis and
  Proposition~\ref{prop:bkwlog}  yield
  ill-posedness properties for the nonlinear Schr\"odinger equation,
  established in \cite{CCT,CCT2}. 
\end{rema}

\subsection{Another instability mechanism}
\label{sec:nonode}

We now consider the general assumptions 
Theorem~\ref{theo:strong} (in
particular, we no longer assume $\phi_0\equiv 0$). 
Seeking $u^h \sim a e^{i\phi /h}$
as $h\to 0$, we now find:
\begin{equation}\label{eq:bkwlim2}
\begin{aligned}
&\d_t \phi +\frac{1}{2}|\nabla \phi|^2 +f\(|a|^2\)=
  0\quad ;\quad 
  \phi(0,x)=\phi_0(x)\,. \\
& \d_t a + \nabla \phi\cdot \nabla a +\frac{1}{2}a \Delta \phi
=0 \quad ; \quad a(0,x)=a_0(x)\, .
\end{aligned}
\end{equation}
As $t\to 0$, approximate $\phi$ and $a$ by their Taylor expansion:
\begin{equation*}
\phi(t,x)\sim \sum_{j\ge 0}t^{j}\phi_j(x)\quad ; \quad a(t,x)\sim
\sum_{j\ge 0}t^{j}{a}_j(x) \, .
\end{equation*}
Now all the powers of $t$ must be taken into account. Using
\eqref{eq:bkwlim2}, we see that $(\phi_{j+1},a_{j+1})$ is given
recursively by $(\phi_l,a_l)_{0\le l \le j}$. Define
\begin{align*}
  \u_k^h(t,x) &=a_0(x) \exp\(\frac{i}{h}\sum_{j=0}^k
  t^{j}\phi_j(x)\),\\
\v_k^h(t,x) &=a_0(x) \exp\(\frac{i}{h}\sum_{j=0}^k
  t^{j}\widetilde\phi_j^h(x)\),
\end{align*}
where $(\widetilde\phi_j^h,\widetilde
a_j^h)_{j\ge 1}$ is constructed like $(\phi_j,a_j)_{j\ge 1}$
with $a_0$ replaced by $\widetilde a_0^h$. By induction, we have: 
\begin{lem}
Under the assumption \eqref{eq:hyppert},
\begin{equation*}
  \|\phi_1 -\widetilde \phi_1^h\|_{H^s}\approx \delta^h\, ,\quad
  \forall s\ge 0\, ,
\end{equation*}
and for any $j\ge 2$ and every $s\ge 0$, 
\begin{equation*}
  \|\phi_j -\widetilde \phi_j^h\|_{H^s}+ \|a_{j-1} -\widetilde
  a_{j-1}^h\|_{H^s}\lesssim \delta^h\, . 
\end{equation*}
\end{lem}\label{lem:1608}
For any $t^h\ll 1$, we have
\begin{equation*}
  \left|u^h(t^h,x) -\u_k^h(t^h,x)\right| \sim |a_0(x)| \left|
\exp \(\frac{i}{h}\phi(t^h,x) - \frac{i}{h}\sum_{j=1}^k \(t^h\)^j \phi_j(x) \)
  -1\right|\, .  
\end{equation*}
This goes to zero provided that $\(t^h\)^{k+1}\ll h$. On the other
hand,
\begin{equation*}
  \left|\u_k^h(t^h,x) -\v_k^h(t^h,x)\right| \sim |a_0(x)| \left|
\exp \(\frac{i}{h}\sum_{j=1}^k \(t^h\)^j \(\phi_j(x)-\widetilde
\phi_j^h(x)\) \)   -1\right|\, .  
\end{equation*}
Since $\widetilde
\phi_0^h=\phi_0$, and from Lemma~\ref{lem:1608}, the main term in the
exponential is
\begin{equation*}
  \frac{t^h}{h}\(\phi_1(x)-\widetilde
\phi_1^h(x)\)\approx \frac{t^h\delta^h}{h}\, .
\end{equation*}
All the other terms are negligible from Lemma~\ref{lem:1608}, 
since $t^h\ll 1$. We then have an instability if:
\begin{equation*}
  \(t^h\)^{k+1}\ll h\quad ;\quad t^h\delta^h\gtrsim h \quad ;\quad
  \(t^h\)^2\delta^h\ll h\, .  
\end{equation*}
All these conditions can be satisfied if we take $k+1 \ge N$. We
conclude this paragraph by showing that in general, the mechanism is
not the same as in the previous section. 

First, if $\phi_0\not \equiv 0$, then trivially $\phi_1$ depends on
$\nabla \phi_0$. If $\phi_0\equiv 0$, we have, for $k\geq 2$, 
\begin{equation*}
  \left| \u_k^h - \u_1^h \right| = |a_0| \left|\exp
    \(\frac{i}{h}\sum_{j=2}^k t^j \phi_j(x)\)-1 \right| 
\end{equation*}
If $\phi_0\equiv 0$, then $\phi_2 \equiv 0$, but in general, $\phi_3
\not \equiv 0$. So if the above instability occurs for $t^h \gtrsim
h^{1/3}$, then it is not an o.d.e. mechanism. We check that if
$\delta^h = h^{2/3}$, then $\u_2^h$ and $\u_1^h$ diverge 
before the instability; therefore, so do $u^h$ and $\u_1^h$.

\subsection{Strong nonlinearities with harmonic potential}
Introduce an isotropic harmonic potential:
\begin{equation*}
ih \d_t u^h + \frac{h^2}{2}\Delta u^h = \frac{|x|^2}{2}u^h +
f\(|u^h|^2\)u^h 
\quad ; \quad u^h(0,x)= a_0(x)e^{i\phi_0(x)/h}\, .
\end{equation*}
Following ideas used in the linear case \cite{Niederer}, we remove the
potential by posing:
\begin{equation}\label{eq:chgtharmo}
U^h(t,x) =
\frac{1}{(1+t^2)^{n/4}}e^{i\frac{t}{1+t^2}\frac{|x|^2}{2h}}u^h \(
\arctan t , \frac{x}{\sqrt{1+t^2}}\)\, .
\end{equation}
Then $U^h$ solves:
\begin{equation}\label{eq:harmomod}\left\{
\begin{aligned}
ih \d_t U^h + \frac{h^2}{2}\Delta U^h &= 
\frac{1}{1+t^2}f\(\(1+t^2\)^{n/2}|U^h|^2\)U^h\, ,\\ 
U^h(0,x)&= a_0(x)e^{i\phi_0(x)/h} .
\end{aligned}\right.
\end{equation}
We can then proceed as above. The only difference is the presence of
time in the nonlinearity, which changes very little at the formal
level. 

\subsection{Weaker nonlinearities}
We come to the framework of Corollary~\ref{cor:weak}: 
\begin{align*}
i\e \d_t u^\e + \frac{\e^2}{2}\Delta u^\e = 
f\(\e^k|u^\e|^2\)u^\e
\quad ; \quad u^\e\big|_{t=0}= a_0(x)e^{-i\frac{|x|^2}{2\e}}\, ,
\end{align*}
where $n\ge 2$, $1<k<n$. Following \cite{Ca-p}, denote $\g =k/n$ and
introduce
\begin{equation}
  \label{eq:pseudo}
  u^\e(t,x) =
  \frac{1}{(1-t)^{n/2}}\psi^\e\(\frac{\e^\g}{1-t}\virgp
  \frac{x}{1-t} \) e^{i\frac{|x|^2}{2\e(t-1)}}\, .
\end{equation}
This can be viewed as a ``semi-classical'' conformal transform, as
compared to the ``usual'' case introduced in \cite{GV82}. 
Then with $h = \e^{1-\g}$, which goes to zero by assumption, and
denoting $t_0^h= h^{\g/(1-\g)}$, $\psi(t,x)$ solves:
\begin{equation}
  \label{eq:psi}
  ih \d_t \psi^h +\frac{h^2}{2}\Delta \psi^h = t^{-2}f\( t^n |\psi^h|^2\)
  \psi^h\quad ; \quad \psi^h\big|_{ t = t_0^h} = a_0(x)\, .
\end{equation}
We can then adapt the preceding approach. This explains
the different notation $\e$ for the semi-classical parameter. Note
that the apparently singular factor $t^{-2}$ is harmless as $t\to 0$,
since we assumed $n\ge 2$ and $f(0)=0$ (this is where this assumption
comes into play). 

Instability
occurs for $\frac{\e^\g}{1-t}\approx t^h$ where $t^h$ and
$\delta^h$ satisfy conditions in the same vein as above. When an
isotropic potential is incorporated, we can essentially 
superimpose the above two changes of unknown functions.

%% file: strong.tex
\section{Proof of Theorem~\ref{theo:strong}}
\label{sec:strong}

\subsection{Case with no potential}
\label{sec:nlsstrong}
In this section, we complete the proof of Theorem~\ref{theo:strong} in
the case $\om=0$. 
Let $n\ge 1$, $a^h_0 \in \S(\R^n)$ bounded in $H^s$ uniformly in
$h\in ]0,1]$ for every $s>0$, and $\phi_0$ as in
Theorem~\ref{theo:strong}.  Consider the initial value problem:
\begin{equation}\label{eq:nlsEmmanuel}
ih \d_t w^h + \frac{h^2}{2}\Delta w^h = 
f\(|w^h|^2\)w^h 
\quad ; \quad w^h(0,x)= a^h_0(x)e^{i\phi_0(x)/h}\, .
\end{equation}
We recall the method of \cite{Grenier98}. It 
somehow boils down to seeking WKB approximation ``the other way
round'': first write the solution as $w^h = \alpha^h e^{i\varphi^h/h}$
(no approximation at this stage), and then study the behavior of
$(\alpha^h,\varphi^h)$ as $h\to 0$, to recover what the usual WKB methods
yield formally. Seek $w^h = \alpha^h e^{i\varphi^h/h}$, with:
\begin{equation}\label{eq:systexact0}
  \begin{aligned}
    \d_t \varphi^h +\frac{1}{2}\left|\nabla \varphi^h\right|^2 + f\(
    |\alpha^h|^2\)= 0\quad &; \quad
    \varphi^h\big|_{t=0}=\phi_0\, ,\\
\d_t \alpha^h +\nabla \varphi^h \cdot \nabla \alpha^h +\frac{1}{2}\alpha^h
\Delta \varphi^h  = i\frac{h}{2}\Delta \alpha^h\quad & ;\quad
\alpha^h\big|_{t=0}= a^h_0\, . 
  \end{aligned}
\end{equation}
Introducing the ``velocity'' ${\tt v}^h = \nabla \varphi^h$,
\eqref{eq:systexact0} yields
\begin{equation}\label{eq:systexact}
  \begin{aligned}
    \d_t {\tt v}^h +{\tt v}^h \cdot \nabla {\tt v}^h + 2 f'\(
    |\alpha^h|^2\) \operatorname{Re}\(\overline{\alpha^h}\nabla \alpha^h\)=
    0\quad &; \quad 
    {\tt v}^h\big|_{t=0}=\nabla \phi_0\, ,\\
\d_t \alpha^h +{\tt v}^h\cdot \nabla \alpha^h +\frac{1}{2}\alpha^h
\operatorname{div}{\tt v}^h  = i\frac{h}{2}\Delta \alpha^h\quad & ;\quad
\alpha^h\big|_{t=0}= a^h_0\, . 
  \end{aligned}
\end{equation}
Separate real and
imaginary parts of $\alpha^h$, $\alpha^h = 
\alpha_1^h + i\alpha_2^h$. Then we have 
\begin{equation}
  \label{eq:systhyp}
  \d_t \bu^h +\sum_{j=1}^n A_j(\bu^h)\d_j \bu^h
  = \frac{h}{2} L 
  \bu^h\, , 
\end{equation}
\begin{equation*}
  \text{with}\quad \bu^h = \left(
    \begin{array}[l]{c}
       \alpha_1^h \\
       \alpha_2^h \\
       {\tt v}^h_1 \\
      \vdots \\
       {\tt v}^h_n
    \end{array}
\right)\quad , \quad L = \left(
  \begin{array}[l]{ccccc}
   0  &-\Delta &0& \dots & 0   \\
   \Delta  & 0 &0& \dots & 0  \\
   0& 0 &&0_{n\times n}& \\
   \end{array}
\right),
\end{equation*}
\begin{equation*}
  \text{and}\quad A(\bu,\xi)=\sum_{j=1}^n A_j(\bu)\xi_j
= \left(
    \begin{array}[l]{ccc}
      {\tt v}\cdot \xi & 0& \frac{\alpha_1 }{2}\,^{t}\xi \\ 
     0 &  {\tt v}\cdot \xi & \frac{\alpha_2}{2}\,^{t}\xi \\ 
     2f'  \alpha_1 \, \xi
     &2f'  \alpha_2\, \xi &  {\tt v}\cdot \xi I_n 
    \end{array}
\right),
\end{equation*}
where $f'$ stands for $f'(|\alpha_1|^2+|\alpha_2|^2)$. 
The matrix $A(\bu,\xi)$ can be symmetrized by 
\begin{equation*}
  S=\left(
    \begin{array}[l]{cc}
     I_2 & 0\\
     0& \frac{1}{4f'}I_n
    \end{array}
\right),
\end{equation*}
which is symmetric and positive since $f'>0$. For an integer
$s>2+n/2$, we bound $(S 
\d_{x}^{\alpha } \bu^h 
, \d_{x}^{\alpha }\bu^h )$ 
where $\alpha $ is a  multi index of length $\le s$, and
$(\cdot ,\cdot)$ is the usual  $L^{2}$ scalar product. We have
\begin{equation*}
\frac{d}{dt}\(S \partial _{x}^{\alpha } \bu^h ,
  \partial _{x}^{\alpha } \bu^h\) 
= \(\partial _{t} S  \partial _{x}^{\alpha } \bu^h , \partial
_{x}^{\alpha } \bu^h\) 
  + 2 \( S \partial _{t}  \partial _{x}^{\alpha } \bu^h , \partial
  _{x}^{\alpha } \bu^h\) 
 \end{equation*}
since $S$ is symmetric. For the first term, we must consider the lower
$n\times n$ block:
\begin{equation*}
\(\partial _{t} S  \partial _{x}^{\alpha } \bu^h ,
  \partial _{x}^{\alpha } \bu^h\) 
\le  \left\|\frac{1}{f'}\d_t\(f'\( | \alpha_1^h|^2 + 
| \alpha_2^h|^2\)\)\right\|_{L^\infty}\(S  \partial _{x}^{\alpha } \bu^h ,
  \partial _{x}^{\alpha } \bu^h\)\, .
\end{equation*}
So long as $\|\bu^h\|_{L^\infty}\le 2\|a_0^h\|_{L^\infty}$, we have:
\begin{equation*}
  f'\( |\alpha_1^h|^2 + 
| \alpha_2^h|^2\) \ge \inf\left\{ f'(y)\ ;\ 0\le y\le
  4\limsup\|a_0^h\|_{L^\infty}^2\right\} =\delta_n>0\, ,
\end{equation*}
where $\delta_n$ is now fixed, since $f'$ is continuous with
$f'>0$. We infer, 
\begin{equation*}
 \left\|\frac{1}{f'}\d_t\(f'\( |\alpha_1^h|^2 + 
| \alpha_2^h|^2\)\)\right\|_{L^\infty}\lesssim \|\bu^h\|_{H^s}\, ,
\end{equation*}
where we used Sobolev embeddings and
\eqref{eq:systhyp}. 
For the second term we use 
\begin{equation*} 
\( S \d_{t}  \d_{x}^{\alpha } \bu^h , \d_{x}^{\alpha } \bu^h\)
= \frac{h}{2}\(S L( \d_{x}^{\alpha } \bu^h )  , \d_{x}^{\alpha }
\bu^h\)
  -  
\Big( S  \d_{x}^{\alpha } \Big(\sum _{j=1}^{n} A_j(\bu^h) \d_{j}
  \bu^h  \Big) , 
 \d_{x}^{\alpha } \bu^h\Big).
\end{equation*} 
We notice that $SL$
is a skew-symmetric second order operator, so the first term is zero. 
For the second term, use the symmetry of $SA_j(\bu^h)$  and usual
estimates on commutators  to get finally:
\begin{equation*}
\frac{d}{dt} \sum _{|\alpha | \le s} \(S \partial _{x}^{\alpha } \bu^h,
\partial _{x}^{\alpha } \bu^h\) \le C\(\left\|\bu^h\right\|_{H^s}\)
\sum _{|\alpha | \le s} \(S \partial _{x}^{\alpha } \bu^h,
\partial _{x}^{\alpha } \bu^h\)\, ,
\end{equation*}
for $s > 2+d/2 $.
Gronwall lemma along with a
continuity argument yield the counterpart of
\cite[Theorem~1.1]{Grenier98}: 
\begin{prop}\label{prop:p}
Let $a^h_0 \in \S(\R^n)$ bounded in $H^m$ uniformly in
$h\in ]0,1]$ for every $m>0$, and let
$s>2+n/2$. Then there exist $T>0$ independent of $h\in ]0,1]$  and
$w^h(t,x) = \alpha^h(t,x) e^{i\varphi^h(t,x)/h}$ solution to
\eqref{eq:nlsEmmanuel} on 
$[0,T]$. Moreover, $\alpha^h$ and
$\varphi^h$ are bounded 
in $L^\infty([0,T];H^s)$,
uniformly in $h\in ]0,1]$. 
\end{prop}
The solution to \eqref{eq:systexact0} formally ``converges'' to the
solution  of:
\begin{equation}\label{eq:systexact00}
  \begin{aligned}
    \d_t \phi^h +\frac{1}{2}\left|\nabla \phi^h\right|^2 + f\(
    |a^h|^2\)= 0\quad &; \quad
    \phi^h\big|_{t=0}=\phi_0\, ,\\
\d_t a^h +\nabla \phi^h \cdot \nabla a^h +\frac{1}{2}a^h
\Delta \phi^h  = 0\quad & ;\quad
a^h\big|_{t=0}= a^h_0\, . 
  \end{aligned}
\end{equation}
The term ``converges'' may not seem appropriate, since the initial
data keeps depending on $h$. Yet, under our assumptions on $a_0^h$,
\eqref{eq:systexact00} has a unique solution $(\alpha^h,\varphi^h)$,
uniformly 
bounded in $L^\infty([0,\tau];H^m)$ for any $m>0$ for some $\tau>0$
independent of $h\in ]0,1]$ (see e.g. \cite{Majda}). We infer:
\begin{prop}\label{prop:estprec}
Let $s\in \N$. There exists $C_s$ independent of $h$ such that for
every $0\le t\le \min(T,\tau)$, 
\begin{equation*}
  \| \alpha^h (t)- a^h(t)\|_{H^s}
  +\| \varphi^h(t) - \phi^h(t)\|_{H^s} \le C_s h t . 
\end{equation*}
\end{prop}
\begin{proof}
  We keep the same notations as above,
  \eqref{eq:systhyp}. Denote by $\bv^h$
  the analog of $\bu^h$ 
  corresponding to $(a^h,\phi^h)$. We have
  \begin{equation*}
\d_t \(\bu^h-\bv^h\) +\sum_{j=1}^n A_j(\bu^h)\d_j \(\bu^h-\bv^h\) +
\sum_{j=1}^n \(A_j(\bu^h)-A_j(\bv^h)\)\d_j \bv^h  
  = \frac{h}{2} L 
  \bu^h    \, .
  \end{equation*}
Keeping the symmetrizer $S$ corresponding to $\bu^h$, we can do
similar computations to  the previous ones. Note that we
know that $\bu^h$ and $\bv^h$ are bounded in
$L^\infty([0,\min(T,\tau)];H^s)$. Denoting ${\bw}^h 
=\bu^h-\bv^h$, we get, for $s>2+n/2$: 
\begin{equation*}
\frac{d}{dt} \sum _{|\alpha | \le s} \(S \partial _{x}^{\alpha } \bw^h,
\partial _{x}^{\alpha } \bw^h\) \lesssim
\sum _{|\alpha | \le s} \(S \partial _{x}^{\alpha } \bw^h,
\partial _{x}^{\alpha } \bw^h\)+ h \| \bw^h(t)\|_{H^s}\, .
\end{equation*}
We conclude with Gronwall lemma.
\end{proof}
This result shows that for small times, WKB solution in the sense of
\eqref{eq:systexact00} provides a good approximation for the exact
solution. Note that since we have to divide phases by $h$, we can
deduce such a result only for times $\ll 1$. The following corollary
is a straightforward consequence of Proposition~\ref{prop:p}:
\begin{cor}\label{cor:1}
Under the assumptions of Proposition~\ref{prop:p}, denote $w_{\rm
  app}^h = a^h e^{i\phi^h/h}$ where $(a^h,\phi^h)$ solves
  \eqref{eq:systexact00}. Then for any $0<t^h\ll 1$, 
\begin{equation*}
\left\| w^h - w_{\rm  app}^h \right\|_{L^\infty([0,t^h];L^2)}\ll 1\,
. 
\end{equation*}
\end{cor}
We now study small time properties of $(a^h,\phi^h)$. 
\begin{defin}\label{def:asym}
  If $T>0$, $(\phi_j^h)_{j\ge 0}$ is a sequence in
  $H^{\infty}(\R^n):=\cap_{s\ge 0}H^s(\R^n)$, and $\phi^h
  \in C([0,T];H^s(\R^n))$ for every $s>0$, the asymptotic 
  relation  
  \begin{equation*}
    \phi^h(t,x)\sim \sum_{j\ge 0}t^j\phi^h_j(x)\quad \text{as }t\to 0
  \end{equation*}
means that for every integer $J\ge 0$ and every $s>0$, 
\begin{equation*}
  \left\| \phi^h(t,\cdot)-\sum_{j= 0}^J t^j\phi^h_j
  \right\|_{H^s(\R^n)} =o\(t^{J}\)\quad \text{as }t\to 0\, .
\end{equation*}
\end{defin}
\begin{prop}\label{prop:small}
  Under the assumptions of Proposition~\ref{prop:p}, there exist
  sequences $(\phi_j^h)_{j\ge 0}$ and $(a_j^h)_{j\ge 1}$ in 
  $H^\infty(\R^n)$ (uniformly in $h\in ]0,1]$), such that the solution of
  \eqref{eq:systexact00} satisfies
  \begin{equation*}
    \phi^h(t,x)\sim  \sum_{j\ge 0}t^j\phi^h_j(x)\ ,\text{ and}\quad
    a^h(t,x)\sim  
    \sum_{j\ge 0}t^{j}a^h_j(x)\quad \text{as }t\to 0\, . 
  \end{equation*}
Moreover, $\phi_0^h=\phi_0$, and $\phi_1^h$ is given by $\phi_1^h
=-f(|a_0^h|^2)$. 
\end{prop}
Plugging such asymptotic series into \eqref{eq:systexact00}, a formal
computation yields a source term which is $\O(t^\infty)$ as $t\to
0$. 
The result then follows with the same approach as in the proof of
Proposition~\ref{prop:p},  
and  Borel lemma (see e.g. \cite{RauchUtah}). Taking
Corollary~\ref{cor:1} into account, we find:
\begin{cor}\label{cor:2}
Let $N\in \N\setminus \{0\}$. Under the assumptions of
Proposition~\ref{prop:p}, for any $0<t^h\ll 
h^{1/(N+1)}$, we have
\begin{equation*}
\left\| w^h - \w_N^h \right\|_{L^\infty([0,t^h];L^2)}\ll 1\, , 
\end{equation*}
where $\w_N^h$ is given by
\begin{equation*}
 \w_N^h(t,x)  = a^h_0(x)\exp \(\frac{i}{h}\sum_{j=0}^N t^j \phi_j^h(x)
 \) .
\end{equation*}
\end{cor}
Applying Corollary~\ref{cor:2} to $u^h$ and $v^h$ respectively yields 
Theorem~\ref{theo:strong} when $\om =0$.

\subsection{With an harmonic potential}
\label{sec:harmostrong}
Now suppose $\om>0$. Up to a
dilation of the coordinates, we can assume that $\om =1$. Let $a^h_0
\in \S(\R^n)$ bounded in $H^s$ uniformly in 
$h\in ]0,1]$ for every $s>0$. Consider the initial value problem:
\begin{equation*}
ih \d_t w^h + \frac{h^2}{2}\Delta w^h = \frac{|x|^2}{2}w^h+
f\(|w^h|^2\)w^h 
\quad ; \quad w^h\big|_{t=0}= a^h_0(x)\, .
\end{equation*}
The change of unknown functions \eqref{eq:chgtharmo} leads to
Equation~\eqref{eq:harmomod} with initial data $a^h_0$. We can then
follow every line of Section~\ref{sec:nlsstrong}. The presence of time in
the nonlinearity does 
not need special care: for the symmetrizer $S$, we can take
\begin{equation*}
  S=\left(
    \begin{array}[l]{cc}
     I_2 & 0\\
     0& \frac{(1+t^2)^{1-n/2}}{4f'}I_n
    \end{array}
\right),\quad \text{where }f'\text{ stands for
}f'\((1+t^2)^{n/2}\(|\alpha_1^h|^2+ |\alpha_2^h|^2\)\).  
\end{equation*}
The presence of time does not perturb the analysis (we always
consider bounded times). We obtain the analogue of
\eqref{eq:systexact00}: 
\begin{equation*}
  \begin{aligned}
    \d_t \phi^h +\frac{1}{2}\left|\nabla \phi^h\right|^2 +
    \frac{1}{1+t^2}f\( \(1+t^2\)^{n/2}
    |a^h|^2\)= 0\quad &; \quad
    \phi^h\big|_{t=0}=0\, ,\\
\d_t a^h +\nabla \phi^h \cdot \nabla a^h +\frac{1}{2}a^h
\Delta \phi^h  = 0\quad & ;\quad
a^h\big|_{t=0}= a^h_0\, . 
  \end{aligned}
\end{equation*}
The conclusions of
Proposition~\ref{prop:small} remain: $\phi_1^h$ is given by the same
formula, but the formulae giving $\(\phi_j^h,a_j^h\)_{j\ge 2}$ 
are different because of time in the nonlinearity.  Since
the change of unknown functions \eqref{eq:chgtharmo} is unitary on
$L^2(\R^n)$, the end of the proof of Theorem~\ref{theo:strong} follows.

%% file: caslim.tex
\section{Proof of Proposition~\ref{prop:caslim}}
\label{sec:caslim}
We study the case with no harmonic potential, $\omega =0$, the case
$\omega>0$ is a straightforward consequence as explained in
Section~\ref{sec:harmostrong}. 

As we noted in Section~\ref{sec:nlsstrong}, the solution to
\eqref{eq:systexact00} yields a good approximation of the solution to
\eqref{eq:nlsEmmanuel} only for \emph{small} times (see
Corollary~\ref{cor:1}). The reason is the same as that mentioned in
Section~\ref{sec:super}: the shift in the cascade of equations in WKB
methods is such that initial corrections of order $h$ become relevant
for times of order $1$. 

For $a_0,a_1 \in\S(\R^n)$ independent of $h$, consider the initial
value problem:
\begin{equation}\label{eq:wcorr}
ih \d_t w^h + \frac{h^2}{2}\Delta w^h = 
f\(|w^h|^2\)w^h 
\quad ; \quad w^h(0,x)= \(a_0(x)+ha_1(x)\)e^{i \phi_0(x)/h}\, .
\end{equation}
We proved in Section~\ref{sec:nlsstrong} that there exists $T_0$
independent of $h\in ]0,1]$ such that $w^h = \alpha^h e^{i\varphi^h/h}$, with
$\alpha^h,\varphi^h \in L^\infty(0,T_0;H^s)$ for every $s\ge 0$, uniformly for
$h\in ]0,1]$. Moreover, Proposition~\ref{prop:estprec} 
yields $(\alpha^h,\varphi^h) =
(a,\phi) + \O(ht)$, where $(a,\phi)$ solves \eqref{eq:bkwlim2}. 
Pursuing the analysis of \cite{Grenier98}, we have:
\begin{prop}\label{prop:correc}
  Let $n\ge 1$, $a_0,a_1 \in\S(\R^n)$,
 $\phi_0 \in C^\infty(\R^n;\R)$ independent of $h$, with $\nabla
 \phi_0 \in H^s(\R^n)$ for every $s\ge 
0$. Let $w^h$ solve \eqref{eq:wcorr}. Define $(a^{(1)},\phi^{(1)})$ by
\begin{equation*}
  \begin{aligned}
    \d_t \phi^{(1)} +\nabla \phi \cdot \nabla \phi^{(1)} +
    2\operatorname{Re}\(\overline a a^{(1)}\)f'\( |a|^2\)&=0\quad
    ;\quad \phi^{(1)}(0,x)=0\, ,\\
   \d_t a^{(1)} +\nabla\phi\cdot \nabla a^{(1)} + \nabla
   \phi^{(1)}\cdot \nabla a + \frac{1}{2} a^{(1)}\Delta \phi
   +\frac{1}{2}a\Delta \phi^{(1)}      &= \frac{i}{2}\Delta a\ ;\
   a^{(1)}(0,x)=a_1 .
  \end{aligned}
\end{equation*}
Then $a^{(1)},\phi^{(1)}\in
L^\infty(0,T_0;H^s)$ for every $s\ge 0$, and
\begin{equation*}
  \|\alpha^h - a - h a^{(1)}\|_{L^\infty(0,T_0;H^s)}+ \|\varphi^h - \phi - h
  \phi^{(1)}\|_{L^\infty(0,T_0;H^s)} \lesssim h^2,\quad \forall s\ge 0\, .
\end{equation*}
The pair $(a,\phi)$ is given by \eqref{eq:bkwlim2}, and does not depend on
$a_1$. 
\end{prop}
The proof is a straightforward consequence of the analysis of
Section~\ref{sec:nlsstrong}, and is given in \cite{Grenier98}. Despite
the notations, it seems unadapted to consider $\phi^{(1)}$ as being
part of the phase. Indeed, we infer from Proposition~\ref{prop:correc}
that 
\begin{equation*}
  \left\|w^h - a e^{i\phi^{(1)}}
    e^{i\phi/h}\right\|_{L^\infty(0,T_0;L^2)}= \O(h). 
\end{equation*}
Relating this information to the WKB methods presented in
Section~\ref{sec:bkw}, we have:
\begin{equation*}
  \a^{(0)} = a e^{i\phi^{(1)}}.
\end{equation*}
Since $\phi^{(1)}$ depends on $a_1$ while $a$ does not, we retrieve
the fact that in super-critical r\'egimes, the leading order amplitude
in WKB methods depends on the initial first corrector $a_1$. Now
Proposition~\ref{prop:caslim} is straightforward, since $(a^{(1)},
\phi^{(1)})$ solves a linear system, and 
\begin{equation*}
  \d_t\phi^{(1)}\big|_{t=0} = -2\operatorname{Re}\(\overline a_0 a_1\)f'\(
    |a_0|^2\).  
\end{equation*}

%% file: weak.tex
\section{Proof of Corollary~\ref{cor:weak}}
\label{sec:weaker}

We indicate how to adapt the analysis of
Section~\ref{sec:strong} when the nonlinearity is attenuated by a
power of the small parameter. By an obvious change of unknown
functions, this 
is equivalent to considering solutions of \eqref{eq:schrodgen2} with
data of order $h^{k/2}$.   
\subsection{Case with no potential}
\label{sec:nlsweak}
Assume $\om=0$. 
For $n\ge 2$, $1<k<n$, $a_0^\e\in\S(\R^n)$ bounded in $H^s$ uniformly in
$\e\in ]0,1]$ for every $s>0$, consider:
\begin{equation*}
i\e \d_t w^\e + \frac{\e^2}{2}\Delta w^\e = 
f\(\e^k|w^\e|^2\)w^\e 
\quad ; \quad w^\e\big|_{t=0}= a_0^\e(x)e^{-i|x|^2/2\e}\, .
\end{equation*}
Introduce $\psi$ given by 
\begin{equation*}
  w^\e(t,x) =
  \frac{1}{(1-t)^{n/2}}\psi^\e\(\frac{\e^\g}{1-t}\virgp
  \frac{x}{1-t} \) e^{i\frac{|x|^2}{2\e(t-1)}}\, .
\end{equation*}
Denoting $\g =k/n$, $h=\e^{1-\g}$ and $t_0^h= h^{\g/(1-\g)}$,
$\psi(t,x)$ solves: 
\begin{equation}
  \label{eq:psi2}
  ih \d_t \psi^h +\frac{h^2}{2}\Delta \psi^h = t^{-2}f\( t^n |\psi^h|^2\)
  \psi^h\quad ; \quad \psi^h\big|_{ t = t_0^h} = a_0^h(x)\, ,
\end{equation}
where we changed the notations $\psi^\e$ and $a_0^\e$ to $\psi^h$ and
$a_0^h$ to keep in mind that these functions depend on the small
parameter. Equation~\eqref{eq:psi2} differs from \eqref{eq:nlsEmmanuel}
by two aspects: the presence of time in the nonlinearity, and the data
are prescribed at time $t=t_0^h$ instead of $t=0$.  

We explain how the computations of Section~\ref{sec:strong} can be
adapted to this case. Seeking $\psi^h =\alpha^h
e^{i\varphi^h/h}$, \eqref{eq:systexact} becomes:
\begin{equation*}
  \begin{aligned}
    \d_t {\tt v}^h +{\tt v}^h \cdot \nabla {\tt v}^h + 2 t^{n-2}f'\(
    t^n |\alpha^h|^2\) \operatorname{Re}\(\overline{\alpha^h}\nabla \alpha^h\)=
    0\quad &; \quad 
    {\tt v}^h\big|_{t=t_0^h}=0\, ,\\
\d_t \alpha^h +{\tt v}^h\cdot \nabla \alpha^h +\frac{1}{2}\alpha^h
\operatorname{div}{\tt v}^h  = i\frac{h}{2}\Delta \alpha^h\quad & ;\quad
\alpha^h\big|_{t=t_0^h}= a^h_0\, . 
  \end{aligned}
\end{equation*}
As a symmetrizer, we take:
\begin{equation*}
  S=\left(
    \begin{array}[l]{cc}
     I_2 & 0\\
     0& \frac{t^{2-n}}{4f'}I_n
    \end{array}
\right),\quad \text{where }f'\text{ stands for }f'\(t^n\(|\alpha_1^h|^2+
|\alpha_2^h|^2\)\). 
\end{equation*}
Unlike in Section~\ref{sec:harmostrong}, we must be careful with the
powers of $t$: the term $t^{2-n}$ in the lower block is singular. 
When computing $\d_t S$ in the energy estimate, differentiating
$t^{2-n}$ on the numerator of the lower block yields a non-positive
term: once again, the assumption $n\ge 2$ is necessary for our proof
to work. When differentiating the denominator, we can factor out 
$(S \d_{x}^{\alpha } \bu^h , \d_{x}^{\alpha }\bu^h )$, times
\begin{equation*}
\left\|\frac{1}{f'}\d_t \( f'\(t^n\(|\alpha_1^h|^2+
|\alpha_2^h|^2\)\)\)\right\|_{L^\infty(\R^n)}.
\end{equation*}
Thus the singular term $t^{2-n}$ is finally harmless.
Apart from that remark, the computations are similar, and we
refer to \cite{Ca-p} for more details. We infer:
\begin{prop}\label{prop:p'}
Let $n\ge 2$, $a^h_0 \in \S(\R^n)$ bounded in $H^m$ uniformly in
$h\in ]0,1]$ for every $m>0$, and let
$s>2+n/2$. Then there exist $T>0$ independent of $h\in ]0,1]$  and
$\psi^h(t,x) = \alpha^h(t,x) e^{i\varphi^h(t,x)/h}$ solution to
\eqref{eq:psi2} on 
$[t_0^h,t_0^h+T]$. Moreover, $\alpha^h$ and
$\varphi^h$ are bounded 
in $L^\infty([t_0^h,t_0^h+T];H^s)$,
uniformly in $h\in ]0,1]$. 
\end{prop}
Similarly, we have the analogue of Corollary~\ref{cor:1}, with:
\begin{equation}\label{eq:bkwweak}
  \begin{aligned}
    \d_t \phi^h +\frac{1}{2}\left|\nabla \phi^h\right|^2 + t^{-2} f\(
    t^n|a^h|^2\)= 0\quad &; \quad
    \phi^h\big|_{t=0}=0\, ,\\
\d_t a^h +\nabla \phi^h \cdot \nabla a^h +\frac{1}{2}a^h
\Delta \phi^h  = 0\quad & ;\quad
a^h\big|_{t=0}= a^h_0\, . 
  \end{aligned}
\end{equation}
Like before, this system has a smooth solution on $[0,\tau]$ for some
$\tau>0$ independent of $h\in ]0,1]$. 
Something must be explained about this approximate system: the time
where data are prescribed is now $t=0$. This seems reasonable since
$t_0^h \to 0$ as $h\to 0$, but there is a price to pay. First, we have
the analogue of Proposition~\ref{prop:small} with different powers of
$t$ due to the presence of time in the nonlinearity, and our
assumption $\phi^h_{\mid t=0}=0$:
\begin{equation}\label{eq:DA}
    \phi^h(t,x)\sim \sum_{j\ge 1}t^{nj-1}\phi^h_j(x)\ ,\text{ and}\quad
    a^h(t,x)\sim  
    \sum_{j\ge 0}t^{nj}a^h_j(x)\quad \text{as }t\to 0\, . 
\end{equation}
To prove the analogue of Corollary~\ref{cor:1}, we compare
$(a^h,\phi^h)\big|_{t=t_0^h}$ with $(a^h,\phi^h)\big|_{t=0}$ thanks to
the above relations. Roughly speaking, the error is of order
$(t_0^h)^{n-1}$.  This yields the following result, whose proof can
be found in \cite{Ca-p}:
\begin{prop}\label{prop:16h35}
Let $n\ge 2$, $a^h_0 \in \S(\R^n)$ bounded in $H^m$ uniformly in
$h\in ]0,1]$ for every $m>0$. 
Let $s\in \N$. There exists $C$ independent of $h$ such that for
every $t_0^h\le t\le \min(T,\tau)$, 
\begin{equation*}
  \| a^h (t)- \alpha^h(t)\|_{H^s}
  +\| \phi^h(t) - \varphi^h(t)\|_{H^s} \le C\( h t +
  h^{\frac{\g(n-1)}{1-\g}}\). 
\end{equation*}
\end{prop}
The last term is $o(h)$ as soon as $k>1$, hence this assumption. Then
we have  the analogue of Corollary~\ref{cor:1}. Using \eqref{eq:DA},
we infer the analogue of Corollary~\ref{cor:2}:
\begin{cor}\label{cor:2'}
Let $N\in \N\setminus\{0\}$. Under the assumptions of
Proposition~\ref{prop:p'}: 
\begin{equation*}
\left\| \psi^h -  \Psi_N^h \right\|_{L^\infty([t_0^h,t^h];L^2)}\ll
1\quad \text{for any }t_0^h\le t^h\ll h^{\frac{1}{(N+1)n-1}}\, , 
\end{equation*}
where $ \Psi_N^h$ is given by:
\begin{equation*}
\Psi_N^h(t,x) = a_0^h(x) \exp\(\frac{i}{h}\sum_{j=1}^N t^{n-1} \phi_j^h(x)\).
\end{equation*}
\end{cor}
We infer Corollary~\ref{cor:weak} in the case $\omega =0$, in the
first case concerning $\delta^\e$. Back to the initial variables, the
instability occurs for 
\begin{equation*}
  \frac{\e^\g}{1-t^\e}\approx \frac{\e^{1-\g}}{\delta^\e}\, ,
\end{equation*}
and the solution at that time is concentrated at scale
$1-t^\e$; $u^\e$ is of order $(1-t^\e)^{-n/2}$. If $\delta^\e
\approx \e^{1-\g -\frac{1}{N}}$, the rate of concentration is then
\begin{equation*}
  1-t^\e \approx \e^{\g -\frac{1}{N}}\, .
\end{equation*}
\subsection{With an harmonic potential}
\label{sec:harmoweak}
With $a_0^\e$ as above, consider now:
\begin{equation*}
i\e \d_t w^\e + \frac{\e^2}{2}\Delta w^\e = \frac{|x|^2}{2}w^\e+
f\(\e^k|w^\e|^2\)w^\e 
\quad ; \quad w^\e\big|_{t=0}= a_0^\e(x)\, .
\end{equation*}
The harmonic potential causes
focusing of the linear solution ($f\equiv 0$) in the limit $\e\to 0$
at time $t=\pi/2$. 
Use the transform \eqref{eq:chgtharmo} to remove the harmonic
potential:
\begin{equation*}
  W^\e(t,x) =
\frac{1}{(1+t^2)^{n/4}}e^{i\frac{t}{1+t^2}\frac{|x|^2}{2h}}w^\e \(
\arctan t , \frac{x}{\sqrt{1+t^2}}\)\, .
\end{equation*}
Now the focusing phenomenon occurs for $W^\e$ when time goes to
infinity. To ``compactify'' time, we use another semi-classical
conformal transform:
\begin{equation*}
\bw^\e(t,x)= \frac{1}{(1-t)^{n/2}}W^\e\(\frac{t}{1-t}\virgp
  \frac{x}{1-t} \) e^{i\frac{|x|^2}{2\e(t-1)}}\, .
\end{equation*}
The focusing for $\bw^\e$ occurs for times close to $1$. It is 
natural to use \eqref{eq:pseudo}:
\begin{equation*}
  \bw^\e(t,x) =
  \frac{1}{(1-t)^{n/2}}\psi^\e\(\frac{\e^\g}{1-t}\virgp
  \frac{x}{1-t} \) e^{i\frac{|x|^2}{2\e(t-1)}},\quad \text{where }\g
  =\frac{k}{n}<1\, .
\end{equation*}
We thus have
\begin{equation*}
  \psi^\e(t,x)=W^\e\(\frac{t}{\e^\g}-1,x \).
\end{equation*}
Keep the notations $h=\e^{1-\g}$ and $t_0^h= h^{\g/(1-\g)}$. The
function $\psi$ solves: 
\begin{equation*}
\left\{
\begin{aligned}
ih\d_t \psi^h +\frac{h^2}{2}\Delta\psi^h &= \frac{1}{(t_0^h)^2 +
  (t-t_0^h)^2 }
f\( \((t_0^h)^2 +  (t-t_0^h)^2\)^{n/2}|\psi^h|^2\)\psi^h,\\
\psi^h\big|_{t=t_0^h} &=a_0^h\, .
\end{aligned}
\right.
\end{equation*}
We have the same equation as \eqref{eq:psi2}, with $t$ in the
nonlinearity replaced by 
$$\((t_0^h)^2 +  (t-t_0^h)^2\)^{1/2}.$$
We can reproduce the analysis of Section~\ref{sec:nlsweak}, with again
\eqref{eq:bkwweak} as a limiting system, since $t_0^h\to 0$ as $h\to
0$. The price to pay is the same: we have an error estimate like in
Proposition~\ref{prop:16h35}, so we must assume $k>1$ to approximate
the phases. 

Back to the initial variables, the
instability occurs for 
\begin{equation*}
  t^\e\approx \frac{\pi}{2}- \arctan\( \frac{\delta^\e}{\e^{1-2\g}}\),
\end{equation*}
and the solution at that time is concentrated at scale
$\cos t^\e$; $u^\e$ is of order $(\cos t^\e)^{-n/2}$. If $\delta^\e
\approx \e^{1-\g -\frac{1}{N}}$, the rate of concentration is then
\begin{equation*}
  \cos t^\e \approx \sin \arctan\( \frac{\delta^\e}{\e^{1-2\g}}\)
      \approx \e^{\g -\frac{1}{N}}\, .
\end{equation*}
In particular, this rate of concentration is large compared to the one
studied in \cite{BZ}, which is $\e$, while the authors consider the
case $n=3$ and $k=2$ (see also Remark~\ref{rem:proj}).

Finally, in the case $\widetilde a_0^\e = a_0 +\e^{1-\g}a_1$,
Corollary~\ref{cor:weak} stems from Proposition~\ref{prop:caslim} in
the same fashion as above. 

%% file: linear.tex
\section{Linear equation}
\label{sec:linear}
In the linear case, justifying WKB methods is
rather easy, and we prove:
\begin{prop}
Let $n\ge 1$, $a_0\in \S(\R^n)$, 
and $V,V_1\in C^\infty(\R^n;\R)$ be  smooth sub-quadratic
potentials:
\begin{equation*}
\d^\alpha V, \d^\alpha V_1\in L^\infty(\R^n)\, ,\quad \forall
\alpha\in\N^n \text{ such that }|\alpha|\ge 2\, .
\end{equation*}
Assume also that $V_1\not \equiv 0$ on $\operatorname{supp}a_0$. 
Let $u^h$ and
$v^h$ solve the initial value problems:
\begin{align*}
ih \d_t u^h + \frac{h^2}{2}\Delta u^h &= V(x)u^h 
\quad ; \quad u^h\big|_{t=0}= a_0(x)\, .\\
ih \d_t v^h + \frac{h^2}{2}\Delta v^h &= \(V(x)+\delta^h V_1(x)\)v^h 
\quad ; \quad v^h\big|_{t=0}= a_0(x)\, .
\end{align*}
Assume that  $h^{2/3}\ll \delta^h \ll
1$. Then we can find $0<t^h\ll h^{1/3}$ such that:
\begin{equation*}
\left\| u^h - v^h \right\|_{L^\infty([0,t^h];L^2)}\gtrsim 1\, .
\end{equation*}
\end{prop}
\begin{proof}
For $\delta \in [0,1]$, let $w^h_\delta$ solve
\begin{align*}
ih \d_t w^h_\delta + \frac{h^2}{2}\Delta w^h_\delta = \(V(x)+\delta
V_1(x)\)w^h_\delta
\quad ; \quad w^h_\delta\big|_{t=0}= a_0(x)\, .
\end{align*}
WKB method yields $w^h_\delta \sim \w^h_\delta = A_\delta
e^{i\Phi_\delta/h}$, where:
\begin{align}
&\d_t \Phi_\delta +\frac{1}{2}|\nabla_x \Phi_\delta|^2 +V(x)+\delta V_1(x)=
  0\quad ;\quad 
  \Phi_\delta(0,x)=0\,. \label{eq:eicon}\\
& \d_t A_\delta + \nabla_x \Phi_\delta\cdot \nabla_x  A_\delta
  +\frac{1}{2}A_\delta  \Delta \Phi_\delta 
=0 \quad ; \quad A_\delta (0,x)=a_0(x)\label{eq:transport}\, .
\end{align}
Since the difference $r^h_\delta:= w^h_\delta - \w^h_\delta$ solves:
\begin{align*}
ih \d_t r^h_\delta + \frac{h^2}{2}\Delta r^h_\delta = \(V(x)+\delta 
V_1(x)\)r^h_\delta
+e^{i\Phi_\delta/h}\frac{h^2}{2}\Delta A_\delta
\quad ; \quad r^h_\delta\big|_{t=0}= 0\, ,
\end{align*}
standard energy estimates yield:
\begin{equation*}
\left\| r^h_\delta\right\|_{L^\infty([0,t];L^2)} \lesssim h \left\|
\Delta A_\delta\right\|_{L^1([0,t];L^2)} \, .
\end{equation*}
Since $V$ and $V_1$ are smooth and sub-quadratic, there exists $T>0$
such that for every $\delta\in [0,1]$, \eqref{eq:eicon} has a smooth
solution $\Phi_\delta\in C^\infty([0,T]\times\R^n)$, and
\eqref{eq:transport} has a smooth solution such that $A_\delta\in
L^\infty([0,T];H^2)$, 
and:
\begin{equation*}
\left\| A_\delta\right\|_{L^\infty([0,T];H^2)}\le C , \text{ where
}C\text{ is independent of }\delta\in [0,1] .
\end{equation*}
Moreover, plugging Taylor expansion in time for $\Phi_\delta$ and
$A_\delta$, we find:
\begin{equation*}
\Phi_\delta(t,x) = -t\(V(x)+\delta V_1(x)\) +\O\(t^3\)\ ;\
A_\delta(t,x) = a_0(x)+\O(t) \text{ as } t\to 0. 
\end{equation*}
This implies that for $0<t^h\ll h^{1/3}$, 
\begin{equation*}
\left\| u^h - v^h\right\|_{L^\infty([0,t^h];L^2)} = \left\| \u^h
- \v^h\right\|_{L^\infty([0,t^h];L^2)} +o(1)\text{ as }h\to 0\, ,
\end{equation*}
where $\u^h$ and $\v^h$ solve the ordinary differential equations:
\begin{align*}
ih \d_t \u^h = V(x)\u^h 
\ ; \ 
ih \d_t \v^h &= \(V(x)+\delta^h V_1(x)\)\v^h 
\ ;\  \u^h\big|_{t=0}= \v^h\big|_{t=0}= a_0(x)\, .
\end{align*}
We infer:
\begin{equation*}
\left\| u^h - v^h\right\|_{L^\infty([0,t^h];L^2)} = \left\| 
a_0(x)\( e^{it\delta^h V_1(x)/h}-1\)\right\|_{L^\infty([0,t^h];L^2)} +o(1)\, .
\end{equation*}
By assumption, we can make the right hand side $\gtrsim 1$ for times
$0<t^h\ll h^{1/3}$ such that 
$t^h\delta^h\gtrsim h$, and the proposition follows. 
\end{proof}

%% file: ill.tex
\section{Application: ill-posedness results}
\label{sec:cct}

As a consequence of the analysis of Section~\ref{sec:ode}, we retrieve
some results established in \cite{CCT,CCT2} concerning ill-posedness
issues for the nonlinear Schr\"odinger equation without a small
parameter. 
\begin{prop}[\cite{CCT2}, \cite{BGTENS}]\label{cor:CCT}
Let $n\ge 1$, $\omega \in \R\setminus \{0\}$ and $\si>0$ an integer.
Consider the nonlinear Schr\"odinger
equation in 
$\R^n$: 
\begin{equation}
  \label{eq:NLSgen}
  i\d_t u +\frac{1}{2}\Delta u = \omega |u|^{2\si}u \quad ;
\quad u_{\mid t=0}=u_0\, .
\end{equation} 
$\bullet$ \emph{Ill-posedness}. Let $s<\frac{n}{2}-\frac{1}{\si}$. Then
\eqref{eq:NLSgen} is not locally well-posed in $H^s(\R^n)$: for any
$\delta >0$, we can find families $(u_{01}^\e)_{0<\e\le 1}$ and
$(u_{02}^\e)_{0<\e\le 1}$  with
\begin{equation*}
  u_{01}^\e,u_{02}^\e\in \S(\R^n)\quad ;\quad \|u_{01}^\e\|_{H^s},
  \|u_{02}^\e\|_{H^s} \le \delta\, ,\quad \|u_{01}^\e -
  u_{02}^\e\|_{H^s}\ll 1\, ,
\end{equation*}
such that if $u_1^\e$ and $u_2^\e$ denote the solutions to
\eqref{eq:NLSgen} with these initial data, there exists $0<t^\e\ll
1$ such that
$\displaystyle
\left\|u_1^\e\(t^\e\)-u_2^\e\(t^\e\) \right\|_{H^s}\gtrsim 1$.\\ 
$\bullet$ \emph{Norm inflation}. Assume  $0<s
<\frac{n}{2}-\frac{1}{\si}$. We can find
$(u^\e)_{0<\e \le 1}$ 
solving 
\eqref{eq:NLSgen}, such that
$u_0^\e \in \S(\R^n)$ and:
\begin{align*}
  \left\|u_0^\e\right\|_{H^s} \ll 1 \quad ; \quad \exists t^\e\ll 1, \
  \left\| u^\e\( 
  t^\e\) \right\|_{H^s} \gg 1\, .
\end{align*}
\end{prop}
\begin{proof}
This result is a 
 straightforward consequence of WKB analysis for small time, as in
 Proposition~\ref{prop:bkwlog}.  For $a_0\in \S(\R^n)$ 
 with $\|a_0\|_{H^s}\le \delta /2$, and $\l >0$, consider $u$ solving
 \eqref{eq:NLSgen} with: 
 \begin{equation*}
   u_0(x) = \l^{-\frac{n}{2}+s}a_0\(\frac{x}{\l}\) .
 \end{equation*}
Using the parabolic scaling and the scaling of $H^s$, define $u_\l$ by:
\begin{equation*}
  u^\l (t,x) = \l^{\frac{n}{2}-s}u\( \l^2 t,\l x\). 
\end{equation*}
It solves:
\begin{equation*}
  i\d_t u^\l +\frac{1}{2}\Delta u^\l = \omega\l^{2-n+2s}|u^\l|^{2\si}u^\l
  \quad ;\quad 
  u_{ \mid t=0}^\l=a_0\,. 
\end{equation*}
Let $h=\l^{\frac{n\si}{2}-1-s\si}$: $\l$ and $h$ go to zero
simultaneously since $s<\frac{n}{2}-\frac{1}{\si}$. Define 
\begin{equation*}
  \psi^h(t,x) = u^\l ( h t,x) =\l^{\frac{n}{2}-s}u\(
  \l^{\frac{n\si}{2}+1-s\si} t,\l x\)\, . 
\end{equation*}
It solves:
\begin{equation}\label{eq:psi44}
  ih\d_t \psi^h +\frac{h^2}{2}\Delta \psi^h =
  \omega|\psi^h|^{2\si}\psi^h \quad 
  ;\quad \psi^h_{\mid t=0} = a_0\, .
\end{equation}
We go back to $u$ via the formula:
\begin{equation*}
  u(t,x) = \l^{-\frac{n}{2}+s}\psi^h\(
  \frac{t}{\l^{\frac{n\si }{2}+1-s\si}},\frac{x}{\l}\) \, .
\end{equation*}
\emph{Ill-posedness}. Let $\widetilde \psi^h$ solve \eqref{eq:psi44}
with a slightly different initial data:
\begin{equation*}
  \widetilde \psi^h_{\mid t=0} = \(1+\delta^h\)a_0\, ,
\end{equation*}
with $\delta^h=|\ln h|^{-\theta}\ll 1$, where $\theta>0$ stems from
Proposition~\ref{prop:bkwlog}. We infer from 
Proposition~\ref{prop:bkwlog} and the discussion of
Section~\ref{sec:ode} that for $t^h=c_0 h|\ln h|^\theta \ll 
1$, we have:
\begin{equation*}
  \left\| \psi^h \( t^h\) - \widetilde\psi^h \( t^h\)\right\|_{H^s}
  \gtrsim 1\, . 
\end{equation*}
Back to the function $u$, this yields the first part of
Proposition~\ref{cor:CCT}. \\

\noindent\emph{Norm inflation}. In \cite{CCT2}, this phenomenon appears as a
transfer of energy from low to high Fourier modes. It corresponds to
the apparition of rapid oscillations in a 
super-critical WKB r\'egime, which can be viewed as a particular case
of the above statement: even though $\psi^h$ is not $h$-oscillatory
initially, rapid oscillations appear instantly. Note that a similar
phenomenon was  shown recently
in the context of Euler equations by C.~Cheverry and O.~Gu\`es
\cite{CG05}. 

Still from Proposition~\ref{prop:bkwlog} and the discussion of
Section~\ref{sec:ode}, with $t^h = c_0 h|\ln h|^\theta$, we
  have
  \begin{equation*}
    \psi^h(t^h)\sim a_0(x)e^{i\frac{t^h}{h}\phi_1(x)}=
    a_0(x)e^{-i\omega|a_0(x)|^{2\si}\(\ln\frac{1}{h}\)^\theta}. 
  \end{equation*}
Even though $\psi^h$ is not yet $h$-oscillatory, ``rapid'' oscillations
have appeared already. Now as in \cite{BGTENS}, we may replace $a_0$
with $ |\ln \l|^{-\theta'}a_0$ to complete the proof of
Proposition~\ref{cor:CCT}. 
\end{proof}

%% file: flow.tex
\section{On the flow map for the cubic, defocusing NLS}
\label{sec:flow}

In the previous section, ill-posedness results were established thanks
to a justification of WKB analysis for very small times, of order
$h |\ln h|^\theta$. For the cubic, defocusing Schr\"odinger
  equation, we saw that a rigorous WKB analysis was available for times of
  order $\O(1)$. 

\begin{proof}[Proof of Corollary~\ref{cor:flow}]
  Mimicking the previous section, for $a_0\in \S(\R^n)$, let
  \begin{equation*}
    u_0(x) = \l^{-\frac{n}{2}+s}a_0\(\frac{x}{\l}\) .
  \end{equation*}
Let $h = \l^{\frac{n}{2}-1 -s}$: $h$ and $\l$ go simultaneously to
zero, since $s<s_c$. Define 
\begin{equation*}
  \psi^h(t,x) = u^\l ( h t,x) =\l^{\frac{n}{2}-s}u\(
  \l^{\frac{n}{2}+1-s} t,\l x\)\, . 
\end{equation*}
It solves:
\begin{equation}\label{eq:psi52}
  ih\d_t \psi^h +\frac{h^2}{2}\Delta \psi^h =
  |\psi^h|^{2}\psi^h \quad 
  ;\quad \psi^h_{\mid t=0} = a_0(x) \, .
\end{equation}
The idea of the proof is that for times of order $\O(1)$, $\psi^h$ is
$h$-oscillatory. This result is expected to be true not only for cubic
defocusing nonlinearities, but it seems this is the only framework where it
has been proved \cite{Grenier98}. 

We infer from Proposition~\ref{prop:correc} that there exist $T>0$
independent of $h\in ]0,1]$, and $a,\phi,\phi_1 \in C([0,T];H^m)$ for
any $m\ge 0$, such that:
\begin{equation*}
  \left\| \psi^h - a e^{i\phi_1} e^{i\phi/ h}\right\|_{L^\infty([0,T];
  H^m)} \le C_m h^{1-m}.
\end{equation*}
Since the $\dot H^m$-norm of $a e^{i\phi_1} e^{i\phi/ h}$ is of order
$h^{-m}$ (when $\phi$ is not stationary), we deduce that there exists
$t\in ]0,T]$ such that for any 
$m\ge 0$:
\begin{equation*}
  \| \psi^h(t)\|_{\dot H^m} \approx h^{-m}.
\end{equation*}
This implies: 
\begin{equation*}
  \left\| u \( \l^{\frac{n}{2}+1 -s}t\)\right\|_{\dot H^k} \approx
  \l^{s-k}\| \psi^h(t)\|_{\dot H^k}\approx \l^{s-k} h^{-k} =
  \l^{s-k-k\(\frac{n}{2}-1 
  -s\)} \, . 
\end{equation*}
The result then follows when considering the limit $\l \to 0$. As in
the previous section, we get exactly the statement of the corollary by
replacing $a_0$ by $|\ln \l|^{-1} a_0$ for instance. 
\end{proof}